\documentclass[12pt]%
{article}
\usepackage{latexsym,amsmath,amssymb}
\usepackage{graphicx,epsf}
\title
{Initial ideals of tangent cones to
	Schubert varieties in orthogonal Grassmannians
	}                    
\author{K.~N.~Raghavan\\
	Institute of Mathematical Sciences\\
        C.~I.~T.~Campus,  Taramani\\
	Chennai 600\,113, INDIA\\
	{\sffamily email\/}: {\ttfamily knr@imsc.res.in}
        \and
        Shyamashree Upadhyay\\
	Chennai Mathematical Institute\\
	Plot No. H1, SIPCOT IT Park\\
	Padur Post, Siruseri 603\,103,	Tamilnadu, INDIA\\
	{\sffamily email\/}: {\ttfamily shyama@cmi.ac.in}}

\usepackage{ifthen}     

\newcommand\wantcolor{no}

\newcommand\finalized{yes}

\newcommand\ignore[1]{}
\usepackage{color}
\providecommand\wantcolor{yes}   %
\ifthenelse{\equal{\wantcolor}{no}}{\renewcommand\color[1]{}}{}
\definecolor{backgroundyellow}{cmyk}{.2,.1,.8,.2}
\definecolor{backgroundblue}{rgb}{0,0,1}
\definecolor{backgroundred}{rgb}{1,0,0}
\definecolor{backgroundmagenta}{cmyk}{0,1,0,0}
\definecolor{GoodForInverseVideo}{rgb}{.6,.8,1}
\definecolor{myyellow}{rgb}{1,1,0}

\newcommand\mysection{\section}
\newcommand\mysectionstar[1]{\section*{#1}}
\newcommand\mysubsection{\subsection}

\newcommand\mysubsubsection[1]{%
		\subsubsection{\sffamily\upshape\mdseries #1}}
\newcommand\mysss{\mysubsubsection}

\newtheorem{annotation}{Annotation}[subsection]
\newtheorem{theorem}[annotation]{
		Theorem}
\newtheorem{lemma}[annotation]{
		Lemma}

\newtheorem{corollary}[annotation]{
		Corollary}
\newtheorem{proposition}[annotation]{
		Proposition}
\newtheorem{example}[annotation]{
		Example}
\newcommand\bexample{\begin{example}\begin{rm}}
\newcommand\eexample{\end{rm}\hfill$\Box$\end{example}}
\newtheorem{notation}[annotation]{
		Notation}
\newcommand\bnotation{\begin{notation}\begin{rm}}
\newcommand\enotation{\end{rm}\end{notation}}

\newtheorem{remark}[annotation]{
		Remark}
\newcommand\bremark{\begin{remark}\begin{sffamily}\begin{upshape}}
\newcommand\eremark{\end{upshape}\end{sffamily}\end{remark}}
\newenvironment{myproof}{%
\par\noindent{\scshape Proof:}\begin{rm}}{\hfill$\Box$\end{rm}\newline}

\numberwithin{equation}{subsection}
\numberwithin{figure}{subsection}

\providecommand\finalized{no}
\ifthenelse{\equal{\finalized}{yes}}{}{\marginparwidth=2cm}
\ifthenelse{\equal{\finalized}{yes}}%
{\newcommand\mylabel[1]{\label{#1}}}%
{\newcommand\mylabel[1]{\label{#1}\marginpar{[{\ttfamily\upshape\tiny #1}]}}}
\usepackage{ragged2e}
\ifthenelse{\equal{\finalized}{yes}}%
{\newcommand\checked[1]{}}%
{\newcommand\checked[1]{\marginpar{[{\ttfamily\upshape\tiny CHECKED: #1}]}}}
\ifthenelse{\equal{\finalized}{yes}}%
{\newcommand\spellchecked[1]{}}%
{\newcommand\spellchecked[1]{\marginpar{[{\ttfamily\upshape\tiny SPELLCHECKED: #1}]}}}

\providecommand\version{public}   
\ifthenelse{\equal{\version}{public}}%
{\newcommand\mcomment[1]{}}%
{\newcommand\mcomment[1]{\marginpar{{\raggedright\sffamily\upshape\tiny #1}}}}
\ifthenelse{\equal{\version}{public}}%
{\newcommand\fcomment[1]{}}%
{\newcommand\fcomment[1]{\footnote{#1}}}

\newcommand\field{\mathfrak{k}}

\newcommand\id{I(d)}

\newcommand\modV{(V)}
\newcommand\misd{\mathfrak{M}_d\modV}
\newcommand\miso{\misd}
\newcommand\sov{{\rm SO}\modV}

\newcommand\pointe{\mathfrak{e}}

\newcommand\yvw{Y(w)}
\newcommand\affinev{\mathbb{A}\modv}

\newcommand\ortho{\ensuremath{\mathfrak{O}}}
\newcommand\ophi{\ensuremath{\ortho\phi}}

\newcommand\pos{\mathfrak{N}}
\newcommand\roots{\mathfrak{R}}

\newcommand\st{\,|\,}

\newcommand\vchain{$v$-chain}
\newcommand\vchains{$v$-chains}

\newcommand\modv{}
\newcommand\posv{\ortho\pos\modv}
\newcommand\oroots{\ortho\roots}
\newcommand\rootsv{\ortho\roots\modv}
\newcommand\andposv{\pos\modv}
\newcommand\androotsv{\roots\modv}

\newcommand\diag{\mathfrak{d}}
\newcommand\diagv{\diag\modv}

\newcommand\mon{\mathfrak{S}}
\newcommand\mont{\mathfrak{T}}

\newcommand\init{\textup{in}_{\torder}}
\newcommand\torder{\vartriangleright}

\newcommand\sign{\textrm{sgn}}
\newcommand\new[1]{\widetilde{#1}}
\newcommand\spnew[1]{\widehat{#1}}

\newcommand\proj[1]{\textup{Proj}\,#1}
\newcommand\projeven[1]{\textup{Proj}^{\textup{e}}\,#1}
\newcommand\card\#

\usepackage{makeidx}
\ifthenelse{\equal{\finalized}{no}}{\makeindex}{}

\ifthenelse{\equal{\finalized}{no}}{
\usepackage[editors,firstabbrev,miniindex]{authorindex}
\usepackage{multicol}
\renewcommand{\cite}{\aicite}    
\aisee{, see }
}{}
\begin{document}
\maketitle
\begin{abstract} 
We compute the initial ideals, with respect to certain conveniently chosen
term orders,  of ideals of tangent cones at torus fixed points
to Schubert varieties in orthogonal Grassmannians.   The initial
ideals turn out to be square-free monomial ideals and therefore 
Stanley-Reisner face rings of simplicial complexes.   We describe
these complexes.   The maximal faces of these complexes encode
certain sets of non-intersecting lattice paths.
\end{abstract}

\vfill\eject
\tableofcontents
\vfill\eject
\mysectionstar{Introduction}\mylabel{s:introduction}
\addcontentsline {toc}{section}{Introduction}
This paper is a sequel to \cite{ru} and the fulfillment of the
hope expressed there that the main result of that paper can be
used to compute initial ideals, with respect to certain `natural'
term orders, of ideals of
tangent cones (at torus fixed points) to Schubert varieties 
in orthogonal Grassmannians.
Any such initial ideal turns out to be generated by square-free monomials
and therefore the Stanley-Reisner face ring of a simplicial complex.
We identify this complex (Theorem~\ref{t:main}). 
The maximal faces of this complex encode a certain set of non-intersecting
lattice paths (Remark~\ref{r:paths}).

The analogous problem for Grassmannians has been addressed
in~\cite{kl1,kr,k,kl2} and for symplectic Grassmannians in~\cite{gr}.
Just as the ideals of tangent cones in those cases are generated
respectively by determinants of generic matrices and determinants of 
generic symmetric matrices,  so the ideals in the present case
are generated by Pfaffians of generic skew
symmetric matrices: see~\S\ref{ss:tangentcone}.
The ideal generated by all Pfaffians
of a fixed degree of a generic skew-symmetric matrix  occurs
as a special case: see~\S\ref{sss:special case}.     
Initial ideals in the special case have been computed in~\cite{ht,jw},
but the term orders there are very different from ours: the Pfaffian
generators are a Gr\"obner basis for those term orders but not for ours.

The present case of orthogonal Grassmannians
features a novel difficulty not encountered 
with either Grassmannians or symplectic 
Grassmannians.   Namely, when one tries,
following the analogy with those cases,  to compute the initial ideal
from the knowledge of the Hilbert function (as obtained in~\cite{ru}),
it becomes evident that, in contrast to those cases,
the natural generators of the ideal of a tangent 
cone---the Pfaffians mentioned above---do {\em not\/} 
form a Gr\"obner basis in any `natural'
term order: see Remark~\ref{r:fail:1}.  
Here what it means for a term order to be `natural'
is dictated by~\cite{ru}:   to each Pfaffian there is naturally
associated a monomial which is a term in it,  and a term order is
{\em natural\/} if the initial term with respect to it of any Pfaffian
is the associated monomial.     This difficulty is overcome
by the main technical result Lemma~\ref{l:main}.

There is another naturally related question that asks if something 
slightly weaker
continues to hold for orthogonal Grassmannians:  namely,
whether the initial ideals of a tangent cone with respect to natural
term orders are all the same.    This too fails: see Remark~\ref{r:fail:2}.
In other words,  the naturalness of a term order turns out not be a
strong determiner,  unlike for ordinary and symplectic Grassmannians.

This paper is organized as follows:
the result is stated 
in~\S\ref{s:theorem} and proved in~\S\ref{s:proof} after
preparations in~\S\ref{s:newforms},~\ref{s:pfaffian}.   There
is heavy reliance on the combinatorial definitions and constructions 
of~\cite{ru}.    Fortunately,  however, only the statement
and not the proof of the main theorem there is used.
\mysection{The theorem}\mylabel{s:theorem}
The whole of this section (except 
for~\S\ref{sss:special case},~\ref{ss:example})
is aimed towards the precise statement of our result,
which appears in \S\ref{ss:result}, after
preparations in~\S\ref{ss:statement}--\ref{ss:termorder}.
For full details about the set up described, see~\cite{ru}.
In~\S\ref{ss:example} the difficulty peculiar to orthogonal Grassmannians
mentioned in the introduction is illustrated by means of an example.
\mysubsection{Initial statement of the problem}\mylabel{ss:statement}
Fix once for all a base field~$\field$ that is algebraically closed 
and of characteristic not equal to~$2$.  
Fix a natural number~$d$, a vector space~$V$ of 
dimension~$2d$, and a non-degenerate symmetric bilinear 
form~$\langle\ ,\ \rangle$ on~$V$.
For $k$ any integer, 
let $k^*:= 2d+1-k$.   Fix a basis 
$e_1$, \ldots, $e_{2d}$ 
of $V$ such that
\[
\langle e_i, e_k \rangle = \left\{ \begin{array}{rl}
    1 & \mbox{ if $i=k^*$ }\\ 
    0 & \mbox{ otherwise}\\ \end{array}\right. \]
Denote by ${\rm SO}(V)$ the group of linear automorphisms of~$V$ that
preserve the form $\langle\ ,\ \rangle$ and also the volume form.
Denote by $\miso'$ the closed sub-variety of the Grassmannian of $d$-dimensional
subspaces consisting of the points corresponding
to isotropic subspaces.     The action of ${\rm SO}(V)$ on~$V$ induces
an action on~$\miso'$.   There are two orbits for this action.  These
orbits are isomorphic:  acting by a linear automorphism that preserves the
form but not the volume form gives an isomorphism.
We denote by $\miso$ the orbit of the span of $e_1$, \ldots, $e_d$ and
call it the {\em (even) orthogonal Grassmannian\/}.   

The {\em Schubert\index{Schubert varieties} varieties\/} of~$\miso$ are defined to be the $B$-orbit
closures in~$\miso$ (with canonical reduced scheme structure), where
$B$ is a Borel subgroup of~$\sov$. 
The problem that is tackled in this paper is this:  given
a point on a Schubert variety in~$\miso$, compute the initial ideal,
with respect to some convenient term order,
of the ideal of functions vanishing on 
the tangent cone to the Schubert variety at
the given point.     The term order is specified in~\S\ref{ss:termorder},
and the answer given in Theorem~\ref{t:main}.

Orthogonal Grassmannians and Schubert varieties in them
can, of course, also be defined when the dimension
of the vector space~$V$ is odd.   As is well known and recalled with
proof in~\cite{ru}, such Schubert varieties are isomorphic to those in
even orthogonal Grassmannians.     The results of this paper would
therefore apply also to them.

\mysubsection{The problem restated}\mylabel{ss:restatement}
We take~$B$ to be the subgroup consisting of elements that are
upper triangular with respect to the basis $e_1,\ldots,e_{2d}$.
The subgroup~$T$ consisting of elements that are diagonal with
respect to $e_1,\ldots,e_{2d}$ is a maximal torus of~${\rm SO}(V)$.
The $B$-orbits of $\miso$ are naturally indexed by its $T$-fixed points:
each orbit contains one and only one such point.     The $T$-fixed points
of~$\miso$ 
are easily seen to be of the form $\langle e_{i_1},\ldots,e_{i_d}\rangle$ for
$\{i_1,\ldots,i_d\}$ in~$I(d)$,  where $I(d)$ is the set of subsets
of $\{1,\ldots,2d\}$ of cardinality~$d$ satisfying the following
two conditions:
\begin{itemize}
\item 
for each $k$, $1\leq k\leq d$, there
does not exist~$j$, $1\leq j\leq d$, such that $i_k^*=i_j$---in other words,
for each $\ell$, $1\leq\ell\leq 2d$,   exactly one
of $\ell$ and $\ell^*$ appears in $\{i_1,\ldots,i_d\}$; 
\item
the parity is even of the number of elements of the subset that are
(strictly) greater than~$d$.
\end{itemize}

Let $I(d,2d)$ denote the set of all subsets of cardinality~$d$ 
of~$\{1,\ldots,2d\}$.    
We use symbols~$v$, $w$, \ldots to denote elements of~$I(d,2d)$ (in particular,
those of~$I(d)$). 
The members of~$v$ are denoted $v_1$, \ldots, $v_d$,
with the convention that $1\leq v_1<\ldots< v_d\leq 2d$.    
There is a natural partial order on~$I(d,2d)$:  $v\leq w$, if
$v_1\leq w_1$, \ldots, $v_d\leq w_d$.   

The point
of the orthogonal Grassmannian~$\miso$ that is the span of
$e_{v_1}$, \ldots, $e_{v_d}$ for $v\in I(d)$ 
is denoted $\pointe^v$.    The $B$-orbit
closure of $\pointe^v$ is denoted~$X(v)$.    The point~$\pointe^v$ (and
therefore the Schubert variety~$X(v)$) is contained in the Schubert
variety~$X(w)$ if and only if~$v\leq w$.

Our problem can now be stated thus:  given elements $v\leq w$ of~$I(d)$,
find the initial ideal of functions vanishing on the tangent cone at
$\pointe^v$ to the Schubert variety~$X(w)$.     The tangent cone being
a subvariety of the tangent space at~$\pointe^v$ to~$\miso$, 
we first
choose a convenient set of co-ordinates for the tangent space.
But for that we need to fix some notation.
\mysubsection{Basic notation}\mylabel{ss:basicnotation}\mylabel{ss:basicnotn}
Let an element~$v$ of~$I(d)$ remain fixed.
We will be dealing extensively with ordered pairs $(r,c)$,
$1\leq r,c\leq 2d$,  such that $r$ is not and $c$ is an entry of~$v$.
Let $\androotsv$\index{rootsv@$\androotsv$} denote the set of all such ordered pairs, and set
\\[1mm]
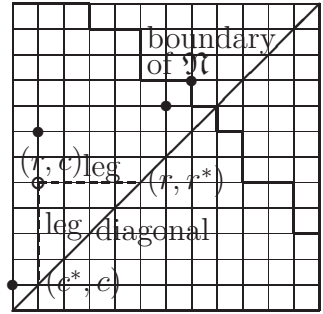
\begin{minipage}{6cm}
\begin{align*}
  \andposv &:= \index{N@$\protect\andposv$}%
  \left\{(r,c)\in\androotsv\st r>c\right\}\\
  \rootsv &:= \index{orootsv@$\protect\rootsv$}%
  \left\{(r,c)\in\androotsv\st r<c^*\right\}\\
  \posv &:= \index{on@$\protect\posv$}%
  \left\{(r,c)\in\androotsv\st r>c, r<c^*\right\}\\
  &=\rootsv\cap\andposv\\
  \diagv &:= \index{diagonal, $\protect\diagv$}%
  \left\{(r,c)\in\androotsv\st r=c^*\right\}\\
\end{align*}\vfill
\end{minipage}
\hfill \begin{minipage}{6cm}
\setlength{\unitlength}{.34cm}
\begin{picture}(12,12)(-3,0)
\multiput(0,0)(12,0){2}{\line(0,1){12}}
\multiput(0,0)(0,12){2}{\line(1,0){12}}
\linethickness{.05mm}
\multiput(0,0)(1,0){13}{\line(0,1){12}}
\multiput(0,0)(0,1){13}{\line(1,0){12}}
\thicklines
\put(0,0){\line(1,1){12}}
\put(3.25,2.75){diagonal}
\thicklines
\put(5.25,10.25){boundary}\put(5.25,9.25){of $\andposv$}
\put(0,12){\line(1,0){3}}
\put(3,12){\line(0,-1){1}}
\put(3,11){\line(1,0){2}}
\put(5,11){\line(0,-1){2}}
\put(5,9){\line(1,0){2}}
\put(7,9){\line(0,-1){1}}
\put(7,8){\line(1,0){1}}
\put(8,8){\line(0,-1){1}}
\put(8,7){\line(1,0){1}}
\put(9,7){\line(0,-1){2}}
\put(9,5){\line(1,0){2}}
\put(11,5){\line(0,-1){2}}
\put(11,3){\line(1,0){1}}
\put(12,3){\line(0,-1){3}}
\put(1,5){\circle{.4}}
\put(.25,5.5){$(r,c)$}
\put(1.25,.75){$(c^*,c)$}
\put(5.25,4.75){$(r,r^*)$}
\put(1.25,3){leg}
\put(2.75,5.25){leg}
\multiput(1,5)(0,-.5){8}{\line(0,-1){.3}}
\multiput(1,5)(.5,0){8}{\line(1,0){.3}}
\put(0,1){\circle*{.4}}
\put(1,7){\circle*{.4}}
\put(6,8){\circle*{.4}}
\put(7,9){\circle*{.4}}
\end{picture}
\end{minipage}
\\

The picture shows a drawing of $\androotsv$.   We think of $r$ and $c$
in $(r,c)$ as row index and column index respectively.    The columns
are indexed from left to right 
by the entries of~$v$ in ascending order,  the rows from top to bottom by
the entries of $\{1,\ldots,2d\}\setminus v$ in ascending order.
The points of $\diagv$ are those on the diagonal, 
the points of $\rootsv$ are those that are (strictly) above the diagonal, and
the points of $\andposv$ are those that are to the South-West of the
poly-line captioned `boundary of $\andposv$'---we draw the
boundary so that points on the boundary belong to $\andposv$.
The reader can readily verify that  $d=13$ and $v=(1,2,3,4,6,7,10,11,13,15,18,19,22)$ for the particular picture drawn.   The points of $\posv$ indicated
by solid circles form a $v$-chain (see~\S\ref{ss:dominate} below).
 
We will be considering {\em monomials}, also called 
{\em multisets\index{multiset@$\textup{multiset}:=\textup{monomial}$}}, in some of these sets.
A {\em monomial}\index{monomial},  as usual,  is a subset with each member being allowed
a multiplicity (taking values in the non-negative integers).  The 
{\em degree\/}\index{degree, of a monomial}
of a monomial has also the usual sense: 
it is the sum of the multiplicities in the
monomial over all elements of the set.  The {\em intersection\/}%
\index{intersection (of a monomial in a set with a subset)}
of a monomial in a set with a subset of the set has also the natural meaning:
it is a monomial in the subset,  the multiplicities being those in the
original monomial.

We will refer to $\diagv\index{diagonal, $\diagv $}$ as the 
{\em diagonal}\index{diagonal, $\diagv $}.   For an element of $\alpha=(r,c)$
of $\androotsv$,  we call $(r,r^*)$ and $(c,c^*)$ its {\em horizontal\/} and
{\em vertical projections\/} (on the diagonal);  they are
denoted by~$p_h(\alpha)$ and $p_v(\alpha)$ respectively.    
For $(r,c)$ in $\posv$,  
its vertical projection belongs to $\andposv$ but not always so
its horizontal projection.
The term {\em projection\/} when not further qualified means either a
vertical or horizontal projection.

\mysubsection{The tangent space to $\miso$ at $\pointe^v$}%
\mylabel{ss:tangentspace}  

Let $\miso\subseteq G_{d}(V)\hookrightarrow \mathbb{P}(\wedge^d
V)$ be the Pl\"ucker embedding (where $G_d(V)$ denotes the
Grassmannian of all $d$-dimensional subspaces of $V$).
For $\theta$ in $I(d,2d)$,  where $I(d,2d)$ denotes the set of subsets
of cardinality~$d$ of $\{1,\ldots,2d\}$, let $p_\theta$ denote 
the corresponding Pl\"ucker coordinate. Consider the affine 
patch $\affinev$ of $\mathbb{P}(\wedge^d V)$
given by $p_v\neq0$, where $v$ is some fixed element of~$I(d)$~($\subseteq
I(d,2d)$).
The affine patch $\affinev^v:=\miso\cap\affinev$ of the orthogonal Grassmannian
$\miso$ is an affine space 
whose coordinate ring can be taken to be the polynomial ring in 
variables of the form $X_{(r,c)}$\index{Xrc@$X_{r,c}$, variable}
with $(r,c)\in\rootsv$.
Taking $d=5$ and $v=(1,3,4,6,9)$ for example,   a general element
of $\affinev^v$ has a basis consisting of column vectors
of a matrix of the following form:
\begin{equation}\label{eq:matrix}\left(\begin{array}{ccccc}
1 & 0 & 0 & 0 & 0 \\
X_{21} & X_{23} & X_{24} & X_{26} & 0 \\
0 & 1 & 0 & 0 & 0 \\
0 & 0 & 1 & 0 & 0 \\
X_{51} & X_{53} & X_{54} & 0 & -X_{26}\\
0 & 0 & 0 & 1 & 0\\
X_{71} & X_{73} & 0 & -X_{54} & -X_{24}\\
X_{81} & 0 & -X_{73} & -X_{53} & -X_{23}\\
0 & 0 & 0 & 0 & 1\\
0 & -X_{81} & -X_{71} & -X_{51} & -X_{21}\\
\end{array}\right)\end{equation}
The origin of the affine space $\affinev^v$, namely the point at which
all $X_{(r,c)}$ vanish,  corresponds clearly to~$\pointe^v$.   The tangent
space to $\miso$ at $\pointe^v$ can therefore be 
identified with the affine space 
$\affinev^v$ with co-ordinate functions $X_{(r,c)}$.
\mysubsection{The ideal~$I$ of 
the tangent cone to $X(w)$ at $\pointe^v$}%
\mylabel{ss:tangentcone}   
Fix elements $v\leq w$ of $I(d)$.
Set $Y(w):=X(w)\cap\affinev^v$,  where $X(w)$ is the Schubert variety 
indexed by~$w$ and $\affinev^v$ is the affine patch around~$\pointe^v$
as in~\S\ref{ss:tangentspace}.
From~\cite{s} we can deduce a set of generators for
the ideal~$I$ of functions on~$\affinev^v$ vanishing on~$Y(w)$
(see for example~\cite[\S3.2.2]{ru}).     We recall this result now.

In the matrix~(\ref{eq:matrix}), columns are numbered by the entries of $v$,
the rows by $1$, \ldots, $2d$.
For $\theta\in\id$, 
consider the submatrix 
given by the rows numbered $\theta\setminus v$ and columns numbered
$v\setminus \theta$.     Such a submatrix being of even size and
skew-symmetric along the anti-diagonal, 
we can define its {\em Pfaffian\/} (see~\S\ref{s:pfaffian}).
Let~$f_\theta$ denote this Pfaffian.    
We have
\begin{equation}\label{eq:ideal}   
I=\left(f_\tau\st \tau\in\id, \tau\not\leq w\right).  \end{equation}

We are interested in the tangent cone to $X(w)$ at~$\pointe^v$ 
or, what is the same,  the tangent cone
to $\yvw\subseteq \affinev^v$ at the origin.
Observe that $f_\theta$
is a homogeneous polynomial of degree the 
$v$-degree of $\theta$, where the {\em $v$-degree}\index{v-degree@$v$-degree}
of $\theta$ is defined as
one half of the cardinality of $v\setminus\theta$.
Because of this, 
$Y(w)$ itself is a cone and so equal to its tangent cone.   
The ideal of the tangent cone 
is therefore the ideal~$I$ in~(\ref{eq:ideal}).
\mysubsubsection{A special case}\mylabel{sss:special case}
The ideal generated by all Pfaffians of a given degree~$r$ of
a generic skew-symmetric $s\times s$ matrix occurs as a special
case of the ideal~$I$ in~(\ref{eq:ideal}):
take $d=s$, $v=(1,\ldots,d)$,  and $w=(2r-1,\ldots,d,2d-2r+3,\ldots,2d)$
($w$ consists of two blocks of consecutive integers).   The initial
ideals in this special case, with respect to certain term orders, have been
computed in~\cite{ht,jw}.    The Pfaffian generators are a Gr\"obner
basis for those orders unlike for ours: see~\S\ref{ss:example}.
\mysubsection{The term order}\mylabel{ss:termorder}
We now specify the term order(s)~$\torder$
on monomials in the co-ordinate functions
(of the tangent space at a torus fixed point) with respect to 
which the initial ideals in our theorem are to be taken.

\newcommand\toone{>_1}
\newcommand\totwo{>_2}
\newcommand\toi{>_i}
Fix an element~$v$ of~$I(d)$.  Let $\toone$ and $\totwo$ be total
orders on~$\rootsv$ satisfying the following conditions.
For both $i=1$ and $i=2$:
\begin{itemize}
\item 
$\alpha\toi\beta$ if $\alpha\in\posv$,  $\beta\in\rootsv\setminus\posv$,
and the row indices of $\alpha$ and $\beta$ are equal;
\item
$\alpha\toi\beta$ if $\alpha\in\posv$,  $\beta\in\posv$,
the row indices of $\alpha$ and $\beta$ are equal,  and the column index
of~$\alpha$ exceeds that of~$\beta$.
\end{itemize}
In addition:
\begin{itemize}
\item 
$\alpha>_1\beta$ (respectively $\alpha<_2\beta$) 
if $\alpha\in\posv$,  $\beta\in\rootsv$ and the row index of $\alpha$
is less than that of~$\beta$.
\end{itemize}
Let~$\torder$ be one of the following term orders on monomials in~$\rootsv$
(terminology as in~\cite[pages~329,~330]{ebud}):
\begin{itemize}
\item the homogeneous lexicographic order with respect to~$\toone$; 
\item the reverse lexicographic order with respect to~$\totwo$.
\end{itemize}
\mysubsubsection{A non-standard possibility for the term order}%
\mylabel{sss:atorder}
Here is another (somewhat non-standard) possibility 
for the term order~$\torder$.
We prescribe it in several steps.   
Let $\mon$ and $\mont$ be distinct monomials in $\oroots$.   
\begin{itemize}
\item 
If $\deg\mon>\deg\mont$ then $\mon\torder\mont$.    
\item
Suppose that $\deg\mon=\deg\mont$.   Then look at
the set of all projections (both vertical and horizontal, 
including multiplicities) on the diagonal
of elements of $\mon$ and $\mont$---some of these
projections may be in $\and\roots\modv$ and not in $\and\pos\modv$.
Let $r_1\geq \ldots\geq r_{2k}$ and $r'_1\geq\ldots\geq r'_{2k}$ 
be respectively the row numbers of these projections
for $\mon$ and $\mont$.    If the two
sequences are different,  then $\mon\torder\mont$ if $r_j>r'_j$ for the least
$j$ such that $r_j\neq r'_j$. 
\item Suppose that the projections on the diagonal of 
$\mon$ and $\mont$ are the same.   Consider the
column numbers of elements in both $\mon$ and $\mont$ 
that give rise to the projection with the least row number (namely
$r_{2k}=r'_{2k}$).    Suppose
$c_1\geq\ldots\geq c_\ell$ and $c'_1\geq\ldots\geq c'_\ell$ are these
numbers respectively for $\mon$ and $\mont$.   If
these sequences are different,  then let $\tilde{j}$ be the least integer $j$ such that $c_j\neq c'_j$. The following three cases can arise:
\begin{itemize}
\item[(a)]
Both $(r_{2k}, c_{\tilde{j}})$ and $(r_{2k}, c'_{\tilde{j}})$ are outside $\posv$.
\item[(b)]
Exactly one of $(r_{2k}, c_{\tilde{j}})$ and $(r_{2k}, c'_{\tilde{j}})$ 
belongs to $\posv$.
\item[(c)]
Both $(r_{2k}, c_{\tilde{j}})$ and $(r_{2k}, c'_{\tilde{j}})$ are inside $\posv$.
\end{itemize}

In case (a), we say that $\mon\torder\mont$ if $c_{\tilde{j}}<c'_{\tilde{j}}$, i.e., $(r_{2k}, c_{\tilde{j}})$ is more towards $\posv$ than $(r_{2k}, c'_{\tilde{j}})$. In case (b), we say that $\mon\torder\mont$ if $(r_{2k}, c_{\tilde{j}})\in\posv$ and $(r_{2k}, c'_{\tilde{j}})\notin\posv$. In case (c), we say that $\mon\torder\mont$ if $c_{\tilde{j}}>c'_{\tilde{j}}$.

If the sequences $c_1\geq\ldots\geq c_\ell$ and $c'_1\geq\ldots\geq c'_\ell$ are the same, then there is an equality of sub-monomials of $\mon$ and $\mont$ consisting of those elements with row numbers $r_{2k}=r'_{2k}$. We remove this sub-monomial from both $\mon$ and $\mont$ and then appeal to an induction on the degree.
\end{itemize}
This finishes the description of the term order~$\torder$.

\mysubsection{$v$-chains and \ortho-domination}\mylabel{ss:dominate}
The description of the initial ideal in
our theorem is in terms of \ortho-domination of monomials.
We now recall this notion from~\cite{ru}.   
%
An element~$v$ of~$I(d)$ remains fixed.

For elements $\alpha=(R,C)$, 
$\beta=(r,c)$  of $\posv$ (or more generally of $\androotsv$),  we write 
$\alpha>\beta$ if $R>r$ and $C<c$.    A sequence $\alpha_1>\ldots>\alpha_k$
of elements of $\posv$ (or of~$\andposv$) 
is called a {\em \vchain\/}.   The points
indicated by solid circles in the picture in~\S\ref{ss:basicnotn} form
a \vchain. 
(For the statement of the theorem we need only consider
\vchains~in $\posv$ but for the proof we will also need \vchains\ in~$\andposv$.
The term `\vchain' without further qualification means
one in~$\posv$.)

To each \vchain~$C$ 
there is associated an element $w_C$ (or $w(C)$)
of~$I(d)$:  see~\cite[\S2.2]{ru}.     An element $w$ of~$I(d)$
{\em \ortho-dominates\/} a \vchain~$C$ if $w\geq w(C)$;  it
{\em \ortho-dominates\/} a monomial $\mon$ in $\rootsv$ 
if it \ortho-dominates every \vchain\ in $\mon\cap\posv$.

\mysubsection{The theorem}\mylabel{ss:result}\mylabel{ss:theorem}
We are now ready to state our theorem.   
Let $\field$ be a field, algebraically closed and of characteristic
not~$2$.
Let $d$ be a positive integer and $\miso$ the (even) orthogonal Grassmannian 
over $\field$~(\S\ref{ss:statement}).
Let $v\leq w$ elements of~$I(d)$,  $X(w)$ the Schubert variety in~$\miso$
corresponding to~$w$, and $\pointe^v$ the torus fixed point in~$\miso$
corresponding to~$v$~(\S\ref{ss:restatement}).
Let $P$ denote the polynomial ring $\field[X_\beta\st\beta\in\oroots]$,
the co-ordinate ring of the tangent space~$\affinev^v$ to~$\miso$ 
at~$\pointe^v$~(\S\ref{ss:basicnotn},~\ref{ss:tangentspace}).
Let~$I$ denote the ideal~(\ref{eq:ideal}) in~$P$ of functions vanishing
on the tangent cone to~$X(w)$ at~$\pointe^v$~(\S\ref{ss:tangentcone}).
Let $\init I$ denote the initial ideal of~$I$ with respect to the term
order~$\torder$~(\S\ref{ss:termorder}).
\begin{theorem}\mylabel{t:main}
  The initial ideal\/ $\init I$ has a vector space basis over~$\field$
consisting of monomials in~$\rootsv$ not \ortho-dominated by~$w$
(\S\ref{ss:dominate}).   In other words,  the quotient ring $P/\init I$
is the Stanley-Reisner face ring of the simplicial complex with
vertices~$\rootsv$ and faces the square-free monomials \ortho-dominated by~$w$.
\end{theorem}
\begin{myproof}
  The main theorem of~\cite{ru} asserts that the dimension as a 
vector space of the graded piece of~$P/I$ of degree~$d$ equals the
cardinality of the monomials in~$\rootsv$ of degree~$d$ that are
\ortho-dominated by~$w$.   Since $P/I$ and $P/\init I$ have the same
Hilbert function,  the same is true with $P/I$ replaced by $P/\init I$.
It is therefore enough to show that every monomial in~$\rootsv$ that
is not \ortho-dominated by~$w$ belongs to~$\init I$, and this 
is proved in~\S\ref{s:lemma}.
\end{myproof}
\bremark\mylabel{r:paths}   The maximal faces of the simplicial
complex, i.e., the square-free monomials in~$\rootsv$ maximal 
with respect to being \ortho-dominated by~$w$, encode a certain
set of non-intersecting lattice paths: see~\cite[Part~IV]{ru}.
\eremark
\mysubsection{An example}\mylabel{ss:example}
Let $v$ in~$I(d)$ be fixed.   To every element $\tau\geq v$ of $I(d)$
there is naturally associated a monomial in~$\posv$ ($\subseteq\rootsv$).   
Namely, with terminology and notation as in~\cite{ru},   it is the result
of the application of the map $\ophi$ to the standard monomial~$\tau$.
This monomial occurs as a term in the Pfaffian~$f_\tau$ defined 
in~\S\ref{ss:tangentcone}.   

\bremark\mylabel{r:fail:1}
Suppose we have a term order $\succ$ on monomials in~$\rootsv$ such
that, for every $\tau\geq v$ in~$I(d)$, 
the initial term of the Pfaffian~$f_\tau$ equals the monomial
associated to~$\tau$ as above:
the term orders~$\torder$ of~\S\ref{ss:termorder} are examples.
 It is natural to expect that,
for $w\geq v$ fixed,  the generators $f_\tau$, $\tau$ in $I(d)$
such that $\tau\not\leq w$, of the ideal~$I$~(\ref{eq:ideal}) form
a Gr\"obner basis with respect to~$\succ$.    The analogous statements
for Grassmannians and symplectic Grassmannians are true~\cite{kr,gr}.
But this expectation fails rather spectacularly (i.e., even in the
simplest examples), as we now observe.
\eremark
   Take $d=5$ and $v=(1,2,3,4,5)$.   Then the top
half of the matrix~(\ref{eq:matrix}) is the identity matrix
and the bottom half looks like this:
\[\left(
  \begin{array}{ccccc}
    a & b & c & d & 0  \\
    e & f & g & 0 & -d \\
    h & i & 0 & -g & -c \\
    j & 0 & -i & -f & -b \\
    0 & -j & -h & -e & -a 
  \end{array}
\right)
\]

Consider the ideal generated by all Pfaffians of degree~$2$ of the
above matrix.   As observed in~\S\ref{sss:special case},  this 
is the ideal~$I$ of~(\ref{eq:ideal}) with $w=(3,4,5,9,10)$.   
There are~$5$ Pfaffians of degree $2$ corresponding to the~$5$ values
of~$\tau$ in $I(d)$ such that $\tau\not\leq w$:  
\[\textrm{$(1,6,7,8,9)$,\quad $(2,6,7,8,10)$,\quad $(3,6,7,9,10)$,\quad $(4,6,8,9,10)$,\quad $(5,7,8,9,10)$.}\]
They are respectively (see Eq.~(\ref{e.pfaffian}))
\[\textrm{$di-cf+bg$,\quad $dh-ce+ag$,\quad 
$dj-be+af$,\quad $cj-bh+ai$,\quad $gj-fh+ei$.}\]
The monomials of $\posv$ attached to the~$5$ elements $\tau$ above 
are respectively
\[\textrm{$di$,\quad $dh$,\quad $dj$,\quad $cj$,\quad $gj$.}\]
The ideal generated by these monomials does not contain any of the
terms in the following element of~$I$:
\begin{equation}\label{e:element}
-h(di-cf+bg)+i(dh-ce+ag)=cfh-bgh-cei+agi.
\end{equation}
So the Pfaffians $f_\tau$ above
are not a Gr\"obner basis with respect
to~$\succ$.

On the other hand,  the initial terms of the Pfaffians~$f_\tau$ above 
with respect to the term order in~\cite{ht} are respectively
\[\textrm{$bg$,\quad $ag$,\quad $af$,\quad $ai$,\quad $ei$}\]
The Pfaffians $f_\tau$ above are a Gr\"obner basis 
with respect to that term order~\cite{ht}. 
\bremark\mylabel{r:fail:2}
The expectation in Remark~\ref{r:fail:1} having failed,  we could ask
whether a weakening of it---also very natural---holds: are
the initial ideals of a tangent cone to~$X(w)$ 
with respect to various natural
term orders all the same (namely, generated by monomials not
\ortho-dominated by~$w$)?     But this too fails as we now observe.
\eremark
Consider the example discussed above.  Identify~$\rootsv=\posv$
with the variables $a$, $b$, \ldots, $j$.    Consider the degree
lexicographic order on monomials in these variables with respect to
a total order on the variables
in which $d$ is bigger than $a$, $b$, $c$, $e$, $f$, $g$;
and $j$ is bigger than $a$, $b$, $e$, $f$, $h$, $i$.    It is readily
verified that this term order is natural in the sense that it
satisfies the condition in Remark~\ref{r:fail:1}:
there are $16$ elements of~$I(d)$: $v$,  the $5$ listed above,
and $10$ others the associated Pfaffians for which are respectively
the $10$ variables.      

Now take a total order that looks like $d>j>a>\ldots$ (the rest can
come in any order).     The corresponding term order picks out
$agi$ as the initial term of the element of~$I'$ in Eq.~(\ref{e:element}),
but the monomial~$agi$ is $\ortho$-dominated by~$w$ as follows readily from
the definitions.
\mysection{New Forms of a $v$-chain}\mylabel{s.newforms}\mylabel{s:newforms}
\newcommand\oldnew{\ddot}
\newcommand\nfone{\oldnew{F}_1}
\newcommand\nftwo{\oldnew{F}_2}
\newcommand\fone{F_1}
\newcommand\ftwo{F_2}
\newcommand\newf{\new{F}}
\newcommand\chainone{{F>D}}
\newcommand\chaintwo{{\oldnew{F}_1>\oldnew{F}_2}}
In this section,  we construct new $v$-chains, called 
{\em new forms\/}, from a given one.
New forms 
play a crucial role in the proof of the main 
Lemma~\ref{l:main}.    In fact, one
may say that their construction, given in~\S\ref{ss:newforms:def} below,
is the main idea in the proof.
A key property of new forms is recorded in~\S\ref{ss:newforms:props}.
In~\S\ref{ss:y} is described 
an association---not that of~\cite{ru}---of an element~$y_C$ 
of~$I(d)$ to a \vchain~$C$.
The elements~$y_C$ also play a crucial in the proof.

An element~$v$ of $I(d)$ remains fixed throughout.
\mysubsection{Some conventions}\mylabel{ss:nf:prepare}
We will often have to compare diagonal elements of $\androotsv$ 
(\S\ref{ss:basicnotn}) with each other.
With regard to such elements, the phrases {\em smaller than\/}
and {\em greater than\/} (and correspondingly the symbols $<$ and $>$)
mean respectively `to the North-East of' and
`to the South-West of'.     We use these phrases in their strict sense
only: `smaller than' means in particular `not equal to'.
This is consistent with the definition of the relation~$>$ on~$\androotsv$
in~\S\ref{ss:dominate}.


With regard to a \vchain\ (whether in $\posv$ or in $\andposv$),  
such terms as `the first element',
`the last element',  `predecessor of a given element' have the
obvious meaning:  in $\alpha_1>\ldots>\alpha_k$, the first element
is $\alpha_1$,  the last $\alpha_k$,  the immediate predecessor
of $\alpha_j$ is $\alpha_{j-1}$, etc.

Two elements $\alpha>\beta$ of $\posv$ are {\em intertwined\/}
if their legs (see the picture in~\S\ref{ss:basicnotn}) intertwine,
or, more precisely,  the vertical projection of~$\beta$ dominates
the horizontal projection of~$\alpha$.    An {\em intertwined
component\/} of a \vchain~$\alpha_1>\ldots>\alpha_m$ 
has the obvious meaning:
it is a block $\alpha_i>\ldots>\alpha_j$ of consecutive elements such that,
$\alpha_k>\alpha_{k+1}$ is intertwined for $i\leq k < j$,   and
$\alpha_{i-1}>\alpha_i$,  $\alpha_j>\alpha_{j+1}$ are not intertwined
(in case $i>1$, $j<m$ respectively).    Clearly a \vchain~$C$ can be
decomposed as $C_1>\ldots>C_\ell$ into its intertwined components.
Observe that, in all intertwined components except perhaps the last,
projections of all elements belong to~$\andposv$.
A \vchain\ is {\em intertwined\/} if it consists of a single intertwined
component.

Let $F$ be an intertwined \vchain.   We define $\proj{F}$ to be
the set (not multiset) of the projections of all its elements on
the diagonal.    Let $\lambda$ be the smallest of all the projections.
Set
\[\projeven{F}:=\left\{ \begin{array}{cl}
    \proj{F} & \quad \textup{if $\proj{F}$ has even cardinality}\\
    \proj{F}\setminus\{\lambda\}
& \quad\textup{otherwise}
    \end{array}\right.\]
For a \vchain~$C$ with intertwined components $C_1>\ldots>C_\ell$,  set
\begin{align*}
\proj{C}&:=\projeven{C_1}\cup\cdots\cup\projeven{C_{\ell-1}}\cup
\proj{C_\ell}\\
\projeven{C}&:= \projeven{C_1}\cup\cdots\cup\projeven{C_{\ell-1}}\cup
\projeven{C_\ell} \end{align*}

For elements $(R,C)$, $(r,c)$ in~$\andposv$,  we say that $(R,C)$
{\em dominates\/} $(r,c)$ if $R\geq r$ and $C\leq c$.  If the elements
belong to the diagonal,  to say $(R,C)$ dominates $(r,c)$ is equivalent
to saying $(R,C)\geq (r,c)$  (see the first paragraph above).   Given \vchains\
$C: \mu_1>\ldots>\mu_m$ and $D:\nu_1>\ldots>\nu_n$ in~$\andposv$,  we say that 
$D$~dominates~$C$ if $n\geq m$ and $\nu_i$ dominates $\mu_i$ for
$i$, $1\leq i\leq m$.

\mysubsection{The construction}\mylabel{ss:newforms:def}\mylabel{ss:nf:const}
Let $E$ be a (non-empty) $v$-chain. 
The construction of a new form depends on two choices.   The first
of these is a {\em cut-off\/}, the choice of an element of~$E$.
Let us write $E$ as $C>D$,  where $C$ is the part of $E$ up-to
and including the cut-off and $D$ the rest of~$E$.
Of course, $D$ can be empty---this happens if and only if the cut-off
is the last element of~$E$---but $C$ is never empty.

Suppose such a cut-off is chosen.
Let us write the $v$-chain $E$ as $C_1>\ldots>C_{\ell-1}>C_\ell>D_1>D_2>\ldots$,
where $C_1>\ldots>C_\ell$ is the decomposition of~$C$ into intertwined
components, $C_\ell>D_1$ is the intertwined component containing~$C_\ell$ 
of $C>D$ (with $D_1$
possibly empty),  and $D_2>\ldots$ is the decomposition of $D\setminus D_1$
into intertwined components.
We will assume in the sequel that $C_\ell$ has at least two elements---one
may also just say that there are no new forms of~$C$ obtained
from the choice of this cut-off in case this condition isn't met.

The {\em new form\/}~$\new{E}$ of $E$ is defined\footnote{
The new form~$\new{E}$ may not always be defined.   
As just remarked, if~$C_\ell$
has only one element then~$\new{C_{\ell}}$ is not defined and so neither
is~$\new{E}$.   As we will see shortly,  $\new{C_{\ell}}$ is not defined
more generally when $\proj{C_{\ell}}$ has evenly many elements and
contains no elements strictly in between the vertical and horizontal
projections of the last element of~$C_\ell$.    
}   to be
$\spnew{C_1}>\ldots>\spnew{C_{\ell-1}}>\new{C_{\ell}}>D_1>\ldots$,
where $\spnew{C_1}$, \ldots, $\spnew{C_{\ell-1}}$, and 
$\new{C_{\ell}}$  
are as described below.
Note that the part $D$ of~$E$ beyond the cut-off 
does not undergo any change.    It will be obvious that (1)~the vertical
projection of the first element does not change in passing from~$C_j$
to $\new{C}_j$ or $\spnew{C_j}$; (2)~the horizontal projection of
the last element gets no smaller in passing from~$C_j$ to $\spnew{C_j}$;
and~(3) the horizontal (respectively vertical) projection of 
the last element gets bigger (respectively no smaller)
in passing from $C_\ell$ to $\new{C_\ell}$.
We are therefore justified in writing~$\new{E}$ as
$\spnew{C_1}>\ldots>\spnew{C_{\ell-1}}>\new{C_{\ell}}>D_1>\ldots$.

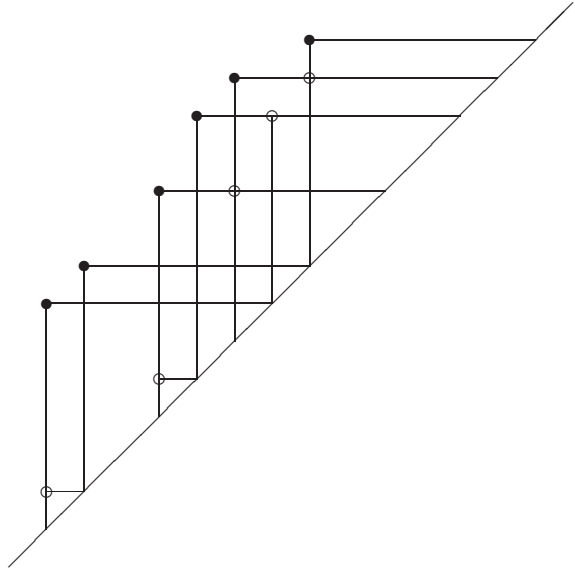
\begin{figure}
\setlength{\unitlength}{.5cm}
\begin{picture}(20,14)(-7,-1)
\linethickness{.2mm}
\put(-1,-1){\line(1,1){15}}
\put(0,1){\circle{.3}}
\put(3,4){\circle{.3}}
\put(5,9){\circle{.3}}
\put(6,11){\circle{.3}}
\put(7,12){\circle{.3}}
\put(0,6){\circle*{.3}}
\put(1,7){\circle*{.3}}
\put(3,9){\circle*{.3}}
\put(4,11){\circle*{.3}}
\put(5,12){\circle*{.3}}
\put(7,13){\circle*{.3}}
\linethickness{.1mm}
\put(0,0){\line(0,1){6}}
\put(1,1){\line(0,1){6}}
\put(3,3){\line(0,1){6}}
\put(4,4){\line(0,1){7}}
\put(5,5){\line(0,1){7}}
\put(6,6){\line(0,1){5}}
\put(7,7){\line(0,1){6}}
\put(0,1){\line(1,0){1}}
\put(0,6){\line(1,0){6}}
\put(3,4){\line(1,0){1}}
\put(1,7){\line(1,0){6}}
\put(3,9){\line(1,0){6}}
\put(4,11){\line(1,0){7}}
\put(5,12){\line(1,0){7}}
\put(7,13){\line(1,0){6}}
\end{picture}
\caption{Illustration of the construction of~$\new{F}$ in the case
when $\proj{F}$ has odd cardinality:  The solid
circles indicate the points of the original \vchain~$F$,  the open
circles those of~$\new{F}$.}\label{f:newodd}
\end{figure}

We first construct $\new{C_\ell}$.  
In fact, we construct $\new{F}$ for an arbitrary intertwined
$v$-chain~$F$ with at least~$2$ elements (subject to a certain
further condition as will be specified shortly).
There are two cases according as the cardinality $\card\proj{F}$
of $\proj{F}$ is odd or even.
Suppose first that it is odd.   In this case no further choice
is involved in the construction. 
Let $(r_1,r_1^*)$, \ldots,$(r_s,r_s^*)$, \ldots, $(r_t,r_t^*)$ be
the elements of~$\projeven{F}$ arranged in decreasing
order,  where $(r_s,r_s^*)$ is
the vertical projection of the last element of~$F$.   
Then~$t$ is even;
and, since there exists at least one horizontal projection that 
is also a vertical projection (because~$\card\proj{F}$ is assumed to be odd),  
we have 
\begin{equation*}\begin{split}
t-s+1 & \leq  \textrm{number of horizontal projections that are}\\
&  \quad\quad\quad
           \textrm{not vertical projections}\\
&<\textrm{number of horizontal projections}\\
&=\textrm{number of vertical projections}\\
& \leq  s
\end{split}
\end{equation*} 
so that $2s-t$ is even and strictly positive.
We define $\new{F}$ to be the~\vchain
\[
(r_2,r_1^*)>\ldots>(r_{2s-t},r_{2s-t-1}^*)>(r_{s+1},r^*_{2s-t+1})>
\ldots>(r_t,r_s^*)\]
In case $s=t$,  the `second half' of~$\new{F}$, namely,
$(r_{s+1},r^*_{2s-t+1})>\ldots>(r_t,r_s^*)$ is understood to be empty.
Figure~\ref{f:newodd} above illustrates the construction.


In the case when $\card\proj{F}$ is even,  the construction of~$\new{F}$
is similar.  The only difference 
is that $(r_1,r_1^*)$, \ldots, $(r_t,r_t^*)$ are now the elements
in decreasing order of the set $\proj{F}$ minus
two elements,  the last element and another that is smaller than
$(r_s,r_s^*)$---if such an element does not exist, then~$\new{F}$
is not defined. 
The choice of such an element 
is the second of the two choices 
involved in the construction of the new form (the first
being the cut-off).  Observe that now~$t-s+2\leq s$,  so that
$2s-t$ is again even and strictly positive.

To define $\spnew{C_1}$,~\ldots,~$\spnew{C_{\ell-1}}$,  we define
more generally $\spnew{F}$ for an arbitrary intertwined \vchain~$F$ both 
projections of all of whose elements belong to~$\andposv$.   Let
$(r_1,r_1^*)$, \ldots, $(r_t,r_t^*)$ be the elements in
decreasing order of $\projeven{F}$.
We define $\spnew{F}$ to be the
\vchain~$(r_2,r_1^*)>\ldots>(r_t,r_{t-1}^*)$.   
\begin{proposition}\mylabel{p:fewer}
With notation as above,  
\begin{enumerate}
\item No two elements of~$\new{C}$ share a projection.
\item $\proj{\new{C}}$ has evenly many elements. It equals\/
$\projeven{C}$ if\/ $\proj{C}$ has oddly many elements.
\item $\new{C}$ has strictly fewer elements than~$C$.
\end{enumerate}
In particular, $\new{E}$ has strictly fewer elements than~$E$.
  \end{proposition}
  \begin{myproof}
    (1) and~(2) being clear from the definition of~$\new{C}$,
we indicate a proof
of~(3).    Using $\card$ to denote cardinality, we have
\[ \card\proj{\new{C}}=
\left\{
  \begin{array}{cl}
\card\projeven{C} & \textrm{if $\card\proj{C}$ is odd}\\
\card\projeven{C}-2 & \textrm{if $\card\proj{C}$ is even}\\
      \end{array}
\right.\]
Because of~(1), $\card\new{C}=\frac{\proj\new{C}}{2}$.   Thus
$\card\new{C}$ equals the greatest integer smaller than
$\frac{\card\proj{C}}{2}$.   But clearly 
$\frac{\card\proj{C}}{2}\leq \card{C}$.
  \end{myproof}
\mysubsubsection{An auxiliary construction}\mylabel{sss:auxiliary}
We now identify a certain sub-\vchain\ of the \vchain~$\new{F}$ 
constructed above.    This auxiliary construction will be used in the proof of
Lemma~\ref{l:p:newdomination:2},  the main ingredient in the proof of
the key property of new forms stated in Proposition~\ref{p:newdomination}.

Let $F>D$ be an intertwined \vchain\ with~
$\new{F}$ being defined.
Let $(r_1,r_1^*)$, \ldots, $(r_s,r_s^*)$, \ldots, $(r_t,r_t^*)$ be as in the
construction of~$\new{F}$ in~\S\ref{ss:newforms:def} above.
Write~$F>D$  
as $F_1>F_2$,  where $F_1$ consists of all elements of~$F$
whose vertical projections 
belong to~$\{(r_1,r_1^*),\ldots,(r_{2s-t},r_{2s-t}^*)\}$
and $F_2$ is the
complement in~$F>D$ 
of~$F_1$.  Denote by~$\oldnew{F}_1$ the part $(r_2,r_1^*)>
\ldots>(r_{2s-t},r_{2s-t-1}^*)$ of~$\new{F}$. 
Consider the
sub-\vchain~$S$ of $\new{F}$ 
consisting of those elements~$(r_j,r_{s-t+j}^*)$, 
$s+1\leq j\leq t$, 
such that
$(r_{s-t+j},r_{s-t+j}^*)$ is the vertical projection of some element 
of~$F_2$ (equivalently of $F_2\setminus D$). 
We set $\oldnew{F}_2$ to be~$S>D$. 
\begin{lemma}\mylabel{l:p:newdom:for}\mylabel{l:aux}
  \begin{enumerate}
\item\mylabel{i:two}
$\oldnew{F}_1>\oldnew{F}_2$ is a sub-\vchain\ of~$\new{F}>D$ 
the inclusion being possibly strict.   
\item\mylabel{i:three}
The projections of $\oldnew{F}_1$ are even in number and all in~$\andposv$.
\item\mylabel{i:four}
The legs of the elements of~$\oldnew{F}_1$ do not intertwine with one another.
Nor does the horizontal leg of the last element of~$\oldnew{F}_1$ intertwine
with the vertical leg of the first element of~$\oldnew{F}_2$.
    \item\mylabel{i:five} The vertical projection of every element of~$F_1$ is a
     projection (vertical or horizontal) of an element of~$\oldnew{F}_1$.
    \item\mylabel{i:six} 
$F_2$ and $\oldnew{F}_2$ are in bijective order preserving correspondence, where
the corresponding elements
      have the same vertical projections (the correspondence is identity
on~$D$).
   Every element of $\oldnew{F}_2$
      has row index no smaller than that of the corresponding element of~$F_2$:
it is bigger for elements of~$\oldnew{F}_2$ not corresponding to elements
of~$D$ (and of course equal for those corresponding to~$D$).

  \end{enumerate}
\end{lemma}
\begin{myproof}
%
(\ref{i:two}) That $\oldnew{F}_1>\oldnew{F}_2$ is a sub-\vchain\
is immediate from the construction.    For an example when it is contained
properly in~$\new{F}$, see Figure~\ref{f:newodd}:  the last but one
open circle does not belong to $\oldnew{F}_1>\oldnew{F}_2$.

(\ref{i:three}) The number of projections of~$\oldnew{F}_1$
is $2s-t$ which is even since $t$~is even.   The horizontal projection
of the last element of~$\oldnew{F}_1$ is $(r_{2s-t},r_{2s-t}^*)$ and this
belongs to~$\andposv$ because $2s-t\leq s$ (since $s\leq t$).

(\ref{i:four}) The first assertion is clear from the definition of~$\oldnew{F}_1$.
The second too is clear: $p_h(\textrm{last element of $\oldnew{F}_1$})
=(r_{2s-t},r_{2s-t}^*) >
(r_{2s-t+1},r_{2s-t+1}^*) 
\geq p_v(\textrm{first element of $\oldnew{F}_2$})$.

(\ref{i:five}) Clear from construction.

(\ref{i:six})  
 Let $F_2$ be $\alpha_1>\ldots>\alpha_k$ and 
$\oldnew{F}_2$ be $\{\beta_1,\ldots,\beta_k\}$,  where $\alpha_i$, $\beta_i$
have the same column index for $1\leq i\leq k$.     
Then $\beta_1>\ldots>\beta_k$,  
for,  $\oldnew{F}_2$ being part of~$\new{F}>D$ 
the $\beta$'s form a \vchain\ 
in some order,  and, their column indices being shared with the $\alpha$'s,
the order $\beta_1>\ldots>\beta_k$ is forced.    

For the second part of the assertion,
let $\alpha_1>\ldots>\alpha_\ell$ be $F_2\setminus D$, and
let $R_1$, \ldots, $R_\ell$
be the respective row indices of $\alpha_1$, \ldots, $\alpha_\ell$.   Then
$r_t> R_\ell$, \ldots, $r_{t-i}> R_{\ell-i}$ for $1\leq i\leq \ell$ 
(for the horizontal projection of the last element of~$F$ 
and possibly one more horizontal
projection  have been discarded from $\proj{F}$ to obtain
$(r_1,r_1^*)$, \ldots, $(r_t,r_t^*)$).    
Also, if~$j$ be such that $(r_j,r_{s-t+j}^*)=\beta_i$ 
for some $i$, $1\leq i\leq \ell$,
then $j\leq t-(\ell-i)$ (strict inequality occurs when  $\new{F}>D$ 
properly contains $\oldnew{F}_1>\oldnew{F}_2$).
    We thus have $r_j\geq r_{t-(\ell-i)}
> R_{\ell-(\ell-i)}=R_i$, which is what we set out to prove.
%
\end{myproof}
\mysubsection{A key property of new forms}\mylabel{ss:newforms:props}
The main result of this subsection is Proposition~\ref{p:newdomination} below.
Invoked in its proof is 
Lemma~\ref{l:p:newdomination:2} which is really where all the action
takes place.

To a \vchain~$C$ of elements in~$\posv$,
there is, as explained in~\cite[\S2.2.2]{ru}, 
an associated element~$w_C$ of~$I(d)$.
There is also a corresponding 
monomial~$\mon_C$ in~$\andposv$ 
associated to~$C$ (\cite[\S5.3.3]{ru}).
%
\bremark\mylabel{r:vchain}  In the statements and proofs of this section
we need to refer to \vchains\ in monomials in~$\andposv$ (typically in
$\mon_C$ where $C$ is a \vchain\ in~$\posv$).   Such \vchains\ are
understood to be in~$\andposv$ (not necessarily restricted to be
in~$\posv$).
\eremark
\begin{proposition}\mylabel{p:newdomination}
 Let $E$ be a \vchain\ in~$\posv$ 
and~$\new{E}$
a new form of~$E$.    Then $w_{\new{E}}\geq w_E$. 
\end{proposition}
\begin{myproof}
By Lemmas~4.5 and~5.5 of~\cite{kr},  it is enough to show that 
every \vchain\ in~$\mon_E$ is dominated by one in~$\mon_{\new{E}}$.
Further, by~\cite[Lemma~5.15]{gr} (or, more precisely, its proof),  it follows,
from the symmetry about the diagonal of monomials attached to \vchains\
in~$\posv$, that it is enough to show that
every \vchain\ in~$\mon_E$ lying (weakly) above the diagonal (in other
words, in~$\posv\cup\diag$) is dominated by one in~$\mon_{\new{E}}$.
We now make some observations after which it will only remain to
invoke Lemmas~\ref{l:p:newdomination:1} and~\ref{l:p:newdomination:2} below.

Decompose $E$ into intertwined components 
$C_1$,~\ldots,~$C_\ell-1$, $C_\ell>D_1$,~\ldots\
as in the description of the construction of the new form~$\new{E}$. 
Let us call these the `parts' of~$E$ (just for now).
There is the corresponding decomposition of~$\new{E}$ into its `parts'
(this is the definition of the {\em parts\/} of~$\new{E}$):
$\spnew{C_1}$, \ldots $\spnew{C_{\ell-1}}$, $\new{C_\ell}>D_1$, $D_2$, \ldots\ .
It is clear from the definitions of $\spnew{C_j}$ and $\new{C_\ell}$ that
each part of~$\new{E}$ is a union of intertwined components.  In particular,
as is immediate from the definition of connectedness in~\S5.3.2 
of~\cite{ru},
each part (of~$E$ or~$\new{E}$) is a union of connected components.  Thus
we have
\[
\mon_E = \mon_{C_1}\cup\cdots\cup\mon_{C_{\ell-1}}\cup\mon_{C_\ell>D_1}\cup
\mon_{D_2}\cup\cdots
\]
and
\[
\mon_{\new{E}} = \mon_{\spnew{C_1}}\cup\cdots\cup\mon_{\spnew{C_{\ell-1}}}\cup\mon_{\new{C_\ell}>D_1}\cup
\mon_{D_2}\cup\cdots\\
\]
Further, since there are no intertwinings between parts,
the following follow easily from the definition of the monomial attached
to a \vchain:
\begin{itemize}
\item 
 any \vchain\ 
$G$ in~$\mon_E$ can be decomposed as: $G_1>\ldots>G_{\ell-1}>G_\ell>H_2>\ldots$
 where $G_1$ is a \vchain\ in~$\mon_{C_1}$, \ldots\ , 
$G_{\ell-1}$ is a \vchain\ in~$\mon_{C_{\ell-1}}$, 
$G_\ell$ is a \vchain\ in~$\mon_{C_\ell>D_1}$,
$H_2$ is a \vchain\ in~$\mon_{D_2}$, \ldots\ ;  
\item
given $v$-chains 
$G_1$ in~$\mon_{\spnew{C_1}}$, \ldots\ , 
$G_{\ell-1}$ in~$\mon_{\spnew{C_{\ell-1}}}$, 
$G_\ell$  in~$\mon_{\new{C_\ell}>D_1}$,
$H_2$ in~$\mon_{D_2}$, \ldots\ ,  all lying weakly above the diagonal,
these can be put together as
$G_1>\ldots>G_{\ell-1}>G_\ell>H_2>\ldots$ to
give a \vchain~$G$ in~$\mon_{\new{E}}$.
\end{itemize}
The proposition now follows from Lemmas~\ref{l:p:newdomination:1} 
and~\ref{l:p:newdomination:2} below.
\end{myproof}
\begin{lemma}\mylabel{l:p:newdomination:1}
For an intertwined \vchain~$F$ both projections of all of whose
elements belong to~$\andposv$,  every \vchain\ in~$\mon_F$ is dominated
by one in~$\mon_{\spnew{F}}$.  (Observe that both~$\mon_F$ and
$\mon_{\spnew{F}}$ consist of diagonal elements.)
\end{lemma}
\begin{myproof}
$\mon_F$ consists of the vertical projections elements of~$F$
in case $\card F$ is even,  and of the vertical projections and the
horizontal projection of the last element in case $\card F$ is odd.
In any case $\mon_F$ consists of evenly many elements.

$\mon_{\spnew{F}}$ consists of all projections of all elements
of~$F$ (in particular, $\mon_{\spnew{F}}\supseteq\mon_F$)
in case the total number of such projections (considered as
a set, not multiset) is even; and, in case that number is odd, 
it consists of all projections except the horizontal projection 
of the last element.  In any case $\mon_{\spnew{F}}$ consists of
evenly many elements.

Suppose that $\mon_{\spnew{F}}\not\supseteq\mon_F$.
Then $\card F$ is odd,   the total number of projections is odd,
and $\mon_F\setminus\mon_{\spnew{F}}=\{\textrm{horizontal projection
of the last element of~$F$}\}$;  in particular, $\card\mon_F=\card F+1$.
Since $\card\mon_{\spnew{F}}\geq\card F$ and  $\card\mon_{\spnew{F}}$ is even, it
follows that $\card\mon_{\spnew{F}}\geq\card F+1$,   
which means that $\mon_{\spnew{F}}$
contains some projection not in~$\mon_F$.     Since any such projection
is bigger than the horizontal projection of the last element of~$F$,  the
lemma follows.
\end{myproof}
\begin{lemma}\mylabel{l:p:newdom:for:2}\mylabel{l:aux:2}
Let $F>D$ be an intertwined \vchain\ with $\new{F}$ being defined.
Let $F_1$, $F_2$, $\oldnew{F}_1$, $\oldnew{F}_2$ be as 
in~\S\ref{sss:auxiliary}.   Then
\begin{enumerate}
\item\mylabel{i:a}
The elements in~$\oldnew{F}_1$ are all of type~H in $\chaintwo$.
\item\mylabel{i:b}
Vertical projections of elements of~$F_1$ belong to~$\mon_\chaintwo$.
\end{enumerate}
\end{lemma}
\begin{myproof}
Statement (\ref{i:a}) follows from (\ref{i:three}) and (\ref{i:four}) 
of~Lemma~\ref{l:aux}.
Statement (\ref{i:b}) from~(\ref{i:a}) and 
Lemma~\ref{l:aux}~(\ref{i:five}).\end{myproof}

\begin{lemma}\mylabel{l:p:newdomination:2}
Let~$F>D$ be an intertwined \vchain\ with~$\new{F}$ being defined.
Given a \vchain\ 
$\mu_1>\mu_2>\ldots$ in $\mon_{F>D}$,  there exists a 
\vchain~$\nu_1>\nu_2>\ldots$ in $\mon_{\new{F}>D}$ that dominates it.
If $\mu_1>\mu_2>\ldots$ lies weakly above the diagonal,  then
$\nu_1>\nu_2>\ldots$ can be chosen also to be so.
\end{lemma}
\begin{myproof}  
Let $F_1$, $F_2$, $\oldnew{F}_1$, $\oldnew{F}_2$ be as defined
in~\S\ref{sss:auxiliary}.  We will show that there exists a 
\vchain~$\nu_1>\nu_2>\ldots$ in~$\mon_{\chaintwo}$ with the desired property.
Since $\oldnew{F}_1>\oldnew{F}_2$
is a sub-\vchain\ of $\newf>D$ (Lemma~\ref{l:aux}~(\ref{i:two})),   this 
will suffice (by either the proof of \cite[Proposition~6.1.1~(1)]{ru} or
\cite[Corollary~6.1.2]{ru} and \cite[Lemmas~4.5,~5.5]{kr}).
For the same reasons as noted in the proof of~Proposition~\ref{p:newdomination},
it is enough to assume that $\mu_1>\mu_2>\ldots$ lies weakly above the
diagonal and find $\nu_1>\nu_2>\ldots$ that dominates it and lies weakly
above the diagonal.    Obviously, we may take without loss of 
generality $\mu_1>\mu_2>\ldots$ to be a maximal such \vchain.

The rest of the proof is divided into three parts:   
\begin{itemize}
\item
Enumerate the maximal \vchains~$\mu_1>\mu_2>\ldots$ in
$\mon_\chainone$ lying weakly above the diagonal.   There are
two of these: see (*) and (**) below.
\item
Identify a certain \vchain~(see~(\dag) below) in~$\mon_\chaintwo$
and lying weakly above the diagonal and list its relevant properties.
\item
Show that the \vchain~(\dag) dominates~(*) in all cases and~(**)
in many cases.   Find a \vchain~(\dag\dag) in $\mon_\chaintwo$ 
and lying weakly above the diagonal that
dominates~(**) when~(\dag) does not.  
\end{itemize}
\newcommand\critical{\sigma}

We start with the first part.
Write $F>D$ as $\alpha_1>\alpha_2>\ldots$ and let~$k$ be the integer
such that $\alpha_k$ is the last element of $F>D$ whose
horizontal projection belongs to $\andposv$: in other words,
$\alpha_k$ is the immediate predecessor of what is called the
critical element in~\cite[\S5.3.4]{ru}.
Of course such an element may not exist, and the proof below, interpreted
properly,  covers that case.   

The \vchain~$F>D$  being intertwined,  its connected components 
(in the sense of~\cite[\S5.3.2]{ru}) are determined by whether or 
not~$\alpha_k$ is connected to its immediate successor:
in either case,  each element $\alpha_j$ for $j\geq k+2$
forms a component by itself,  and the elements $\alpha_1$, \ldots, $\alpha_k$
are all in a single component.   Consider the types of elements 
of~$F>D$ as in~\cite[\S5.3.4]{ru}.
The possibilities for the sequence of these are listed in the following
display.  In these, the underlined type is that of the 
element~$\alpha_k$, the overlined type is that of 
either~$\alpha_k$ or its immediate 
predecessor~$\alpha_{k-1}$ according as whether $k$ is odd or even,
and the vertical bar indicates where the first disconnection
occurs (either just after $\alpha_k$ or just after $\alpha_{k+1}$): 
\newcommand\typev{\textup{V}}
\newcommand\typeh{\textup{H}}
\newcommand\types{\textup{S}}
\newcommand\ul\underline
\newcommand\ol{\overline}
\newcommand\disconnect{\st}
\[
\begin{array}{clllllllll}
\textrm{Case I:} &\ &\typev &\ldots &\typev\,\, &\ul{\ol{\typeh}}\st &{\types} &\types &\types &\ldots \\[2mm]
\textrm{Case II:} &\ &\typev &\ldots &\typev &\ul{\ol{\typev}} &{\typev}\st &\types &\types &\ldots\label{case2} \\[2mm]
\textrm{Case III:} &\ &\typev &\ldots &\typev &{\ol{\typev}} &\ul{\typev}\st &\types &\types &\ldots \\[2mm]
\textrm{Case IV:} &\ &\typev &\ldots &\typev &{\ol{\typev}} &\ul{\typev} &\types\st &\types &\ldots \\[2mm]
\end{array}\]
That these possibilities are all follows readily from the definition of type.

For an element~$\lambda$ of a \vchain~$C$ (in~$\posv$), let
$q_{C,\lambda}$ denote $p_v(\lambda)$ if $\lambda$ is of type~V or~H
and $\lambda$ itself if it is of type~S.   It is easy to see
(and in any case explicitly stated in~\cite[Proposition~5.3.4~(1)]{ru})
that
$q_{C,\lambda}>q_{C,\lambda'}$ for (not necessarily consecutive) 
elements $\lambda>\lambda'$ in~$C$.
It follows that, in Cases~II,~III, and~IV,
\[ (*)\quad\quad q_{\chainone,\alpha_1}>q_{\chainone,\alpha_2}>\ldots\]
is the unique maximal \vchain\ in $\mon_\chainone$ lying weakly
above the diagonal;
in Case~I too it is a maximal \vchain\ but there is
also another one, namely, \[(**)\quad\quad 
p_v(\alpha_1)>p_v(\alpha_2)>\ldots>p_v(\alpha_{k})
>p_h(\alpha_{k})\]
(if $p_h(\alpha_{k})$ dominated $\alpha_j$ for
some $j$, $k < j$, it would contradict the disconnection
between $\alpha_{k}$ and $\alpha_{k+1}$: recall that $\alpha_{k}$ and $\alpha_{k+1}$ are intertwined).   This finishes our first
task of determining the maximal \vchains\ in~$\mon_{\chainone}$ that
lie weakly above the diagonal.

Next we identify a certain \vchain~(see (\dag) below) in $\mon_\chaintwo$
that will have the desired property in almost all cases.
Let $e$ be the integer such that $F_1$ is $\alpha_1>\ldots>\alpha_e$ 
(and $F_2$ is $\alpha_{e+1}>\ldots$). 
Let $\beta_{e+1}>\ldots$ be the counterparts in $\oldnew{F}_2$   respectively of
$\alpha_{e+1}>\ldots$,   the correspondence $\alpha\leftrightarrow\beta$
being as in~Lemma~\ref{l:aux}~(\ref{i:six}):      
\begin{enumerate}
\item[(a)] 
The vertical projections of $\alpha_j$ and $\beta_j$ are
equal for $j= e+1, e+2, \ldots$.   And the row index of~$\beta_j$ is
no less than that of $\alpha_j$ (Lemma~\ref{l:aux}~(\ref{i:six})).
\end{enumerate}

Let $f$ be the largest integer, $f\geq e$, such that $\beta_f$
is of type~V or~H in $\chaintwo$: if either $\alpha_{e+1}$
does not exist or $\beta_{e+1}$ is of type~S,  then $f:=e$ and
$\beta_e$ is taken to be the last element of $\oldnew{F}_1$
(this is not to say that the cardinality of~$\oldnew{F}_1$ is~$e$).
%
Consider the subset~$Z$ of~$\mon_\chaintwo$ consisting of contributions
of elements up to and including~$\beta_f$ 
and only those contributions that are not smaller than~$p_v(\beta_{f+1})$ 
(equivalently~$\beta_{f+1}$): if $\beta_{f+1}$ does not exist, then this
condition is vacuous.   In other words,  $Z$ consists of (1)~the
vertical projections of all elements of $\chaintwo$ up to and 
including~$\beta_f$;  and (2)~the horizontal projections of all elements
of~$\chaintwo$ of type~H except perhaps of $\beta_f$ itself:   the horizontal
projection of $\beta_f$ does not belong to~$Z$ if it is smaller than
$p_v(\beta_{f+1})$ (even if~$\beta_f$ should be of type~H).  Letting 
the elements of~$Z$ arranged in order be $\gamma_1>\ldots>\gamma_g$,
we have the following \vchain\ in~$\mon_\chaintwo$:
\[ (\dag)\quad\quad \gamma_1>\ldots>\gamma_g>\beta_{f+1}>\beta_{f+2}>\ldots \]

We claim:
\begin{enumerate}
\item[(i)]
$p_v(\alpha_1)$, \ldots, $p_v(\alpha_f)$ belong to~$Z$. (So $g\geq f$.)
\item[(ii)]
The horizontal projection of $\alpha_{f+1}$  
does not belong to~$\andposv$.   That is, $f\geq k$ with~$k$ as defined
earlier.
\item[(iii)]
The types of $\alpha_{f+2}$, $\alpha_{f+3}$, \ldots\ in~$\chainone$ are all~S.  
\item[(iv)]
The type of $\alpha_{f+1}$ in~$\chainone$ is either~V or~S.  If it is~V,
then $f=k$ and we are in Case~II (in the enumeration of types listed above).
\item[(v)]
The critical element of~$\chaintwo$ (if it exists) is either~$\beta_f$
or~$\beta_{f+1}$.
\item[(vi)]
If $g\not\geq f+1$ (observe that $g\geq f$ always by (i)),
then $e$ is even. 
\item[(vii)] If $g\not\geq f+1$ and $f$ is odd, 
then~$\beta_f$ is of type~H (in~$\chaintwo$) 
and $\alpha_{f+1}$ is of type~S (in $\chainone$, if $\alpha_{f+1}$ exists).
\end{enumerate}
\begin{myproof}
(i)~If $j\leq e$ (i.e., if $\alpha_j$ belongs to~$F_1$),  
then $p_v(\alpha_j)$ belongs to $Z$ by Lemma~\ref{l:aux:2}~(2);
if $e<j\leq f$,  then $p_v(\alpha_j)=p_v(\beta_j)$ (see~(a) above) 
and so belongs to~$Z$.

(ii)~On the one hand, $p_h(\beta_{f+1})\not\in\andposv$,
for $\beta_{f+1}$ is of type~S.
On the other hand, the row index of 
$\beta_{f+1}$ is at least that of~$\alpha_{f+1}$
(see~(a) above).

(iii) and (iv) follow from combining (ii) with the enumeration of cases
of types of elements of~$\chainone$ above (Cases~I--IV).

(v)~This follows from the definition of type and the choice of~$f$:
an element of type~S cannot precede the critical element;
an element of type~V cannot succeed the critical element.

(vi) 
Suppose that $e$~is odd.
The contributions to $\mon_\chaintwo$ of
elements of~$\oldnew{F}_1$ include $p_v(\alpha_1)$,
\ldots, $p_v(\alpha_e)$ and are evenly many in number 
(Lemma~\ref{l:aux:2}~(1)); $Z$ 
contains all of these (Lemma~\ref{l:aux}~(\ref{i:four})) in addition to
$p_v(\beta_{e+1})$, \ldots, $p_v(\beta_f)$,  so $g\geq (e+1) + (f-e)
=f+1$.   Thus~$e$ is even.

(vii)~By~(vi),  $e$ is even.   Since $f$ is odd, it follows that
$f\geq e+1$.    We first show that $h$ is odd,  where $\beta_h$ is the
first element of the connected component of~$\chaintwo$ that contains~$\beta_f$.
Consider a connected component of 
$\oldnew{F}_2>D$ contained entirely within $\{\beta_{e+1},\ldots,\beta_{f-1}\}$
(if any should exist) (if $f=e+1$,  then $\{\beta_{e+1},\ldots,\beta_{f-1}\}$
is understood to be empty).   If its cardinality is odd, then its last
element, say~$\beta_i$, 
has type~H (this follows from the definition of type: by choice 
of~$f$,  the type can only be~V or~H),   and~$p_h(\beta_i)$ is bigger 
than~$p_v(\beta_{i+1})$
(for otherwise~$\beta_{i+1}$ will be forced to have type~S 
(\cite[Proposition~5.3.4~(1) and~(3)]{ru}),  a contradiction to
the definition of~$f$); and~$Z$ would contain
$p_h(\beta_i)$ in addition to the elements in~(i),  a contradiction.
Thus all such components have even cardinality.    This implies that
$h-e$ is odd, and, since $e$ is even (by (vi)), that~$h$ is
odd.

Since $\beta_{f+1}$ is of type~S (by choice of~$f$),  it is the
last element in its connected component and the component has odd
cardinality.   Since $h$ and~$f$ are odd,  this component can only
be~$\{\beta_{f+1}\}$.    This means that $\beta_f$ is the last element
in its connected component,  and so of type~H:  its type is either~V or~H
by choice of~$f$,  and further because $f-h+1$ is odd its type is~H.

If $p_h(\beta_f)\geq p_v(\beta_{f+1})$,  then $g\geq f+1$,  for $Z$ would
contain $p_h(\beta_f)$ in addition to the elements in~(i).  
So $p_h(\beta_f)<p_v(\beta_{f+1})$. Since $\beta_{f+1}$ is not
connected to~$\beta_f$ (as was just shown), it follows that
$R'\leq R^*$ where $R$, $R'$ are the row indices of~$\beta_f$,~$\beta_{f+1}$.
Letting $r$, $r'$ be the row indices of $\alpha_f$,~$\alpha_{f+1}$,  we
have, by (a)~above, $r'\leq R'\leq R^*\leq r^*$.
This means that $\alpha_{f+1}$ is not connected to~$\alpha_f$ and so
is of type~S (see~(ii) above).
\end{myproof}

The second part of the proof (of the lemma) being over, we start on the third.
We first show that~(\dag) dominates~(*).    
From~(a) above and (iii)~of the claim, it follows that
$q_{\chaintwo,\beta_{f+2}}=\beta_{f+2}>
q_{\chaintwo,\beta_{f+3}}=\beta_{f+3}> \ldots $ dominates
$q_{\chaintwo,\alpha_{f+2}}=\alpha_{f+2}>
q_{\chaintwo,\alpha_{f+3}}=\alpha_{f+3}>\ldots$.
From~(i) of the claim it follows that $\gamma_1>\ldots>\gamma_g>
q_{\chaintwo,\beta_{f+1}}$ dominates $q_{\chainone,\alpha_1}>\ldots>
q_{\chainone,\alpha_{f+1}}$ if either $q_{\chaintwo,\beta_{f+1}}$ dominates
$q_{\chainone,\alpha_{f+1}}$ (which fails by~(a) only when $\alpha_{f+1}$ 
has type~V)
or $g\geq f+1$ (by the definition of~$Z$ and~(a)).
Suppose that $\alpha_{f+1}$ has type~V.
It follows from~(iv) of the claim that $f$ is odd,
and so, from~(vii) of the claim, that $g\geq f+1$.
Thus (\dag) dominates (*).

Now assume that the types of the elements of $F>D$ are as in Case~I
and that $\mu_1>\mu_2>\ldots$ is~(**).   If $f\geq k+1$,   then (\dag) 
dominates (**),  for (\dag) contains
$p_v(\alpha_1)$, \ldots, $p_v(\alpha_k)$, $p_v(\alpha_{k+1})$ (see~(i) 
of the claim),  and $p_v(\alpha_{k+1})\geq p_h(\alpha_k)$
(for $F>D$ is intertwined);   so assume that $f=k$ (by~(ii), we have
$f\geq k$ always).
If $g\geq f+1=k+1$,  then again (\dag) dominates~(**) for similar
reasons: $Z$ contains $p_v(\alpha_1)$, \ldots, $p_v(\alpha_k)$, and it
also contains $g$ elements that 
dominate~$p_h(\alpha_k)$: $p_v(\beta_{k+1})=
p_v(\alpha_{k+1})\geq p_h(\alpha_k)$ for $F>D$ is intertwined.
So assume that $g=f=k$ ($g\geq f$ always by (i)).   
Since we are in Case~I,  $k$ is odd (and hence so is~$f$).    
By~(vii), $\beta_f$ is of type~H
and the following \vchain\ is in~$\mon_\chaintwo$:
\begin{multline*}\textup{(\dag\dag)}\quad 
p_v(\alpha_1)>\ldots>p_v(\alpha_e)> 
p_v(\alpha_{e+1})(=p_v(\beta_{e+1}))> \\
              \ldots> p_v(\alpha_f)(=p_v(\beta_f))>p_h(\beta_f)
\end{multline*}
This \vchain\ dominates (**) by
(a) above.
\end{myproof}

\mysubsection{The element~$y_E$ attached to a \vchain~$E$}%
\mylabel{ss:y}
Let~$E$ be a \vchain\ in~$\posv$.   From~$\projeven{E}$ we can get
an element~$y_E$ of~$I(d,2d)$ by the following natural process
(see the proof of~\cite[Proposition~4.3]{kr}):
the column indices of elements of~$\projeven{E}$
occur as members of~$v$;   these are replaced by the row indices to
obtain~$y_E$.  
\begin{proposition}\mylabel{p:ygeqv}
$y_E\geq v$ and $y_E$ belongs to $I(d)$.
\end{proposition}
\begin{myproof}
Think of~$y_E$ as being the result of a series of operations done starting
with~$v$.    
Let $x\in I(d)$ be
such that $x\geq v$.   Suppose $(r,c)\in\posv$ is such that $c$ occurs
and~$r$ does not in~$x$.
Let $x'$ be the result of replacing~$c$ and~$r^*$ in~$x$ by~$r$ and~$c^*$.   
Then, clearly, either $r>r^*$ in which case $r^*\leq d<d+1\leq c^*$
and $c\leq d<d+1\leq r$,  or $r<r^*$ in which case $c<r\leq d<d+1\leq r^*<c^*$.
In either case $x'\geq x\geq v$ and $x'$ belongs to~$I(d)$.

The proposition follows easily, as we now show,
from the observation just made.     Consider the elements of~$\projeven{E}$
that are not in~$\andposv$.    These can only be horizontal projections,
each of some unique element of~$E$.   Pair these up,  each with the
vertical projection of the corresponding element of~$E$
(all vertical projections belong to~$\projeven{E}$).  
Since $\projeven{E}$ has even cardinality,  there are evenly
many elements left (all in~$\andposv$) after the elements not in~$\andposv$
are paired up as prescribed.  Pair these up in some arbitrary way.     
If $(r,r^*)$
and $(c^*,c)$ are the horizontal and vertical projections of an
element~$(r,c)$ in~$\posv$,  we can think of replacing $r^*$ by $r$ and
$c$ by $c^*$ as the single operation described in the previous paragraph
in going from~$x$ to~$x'$.
It should now be clear that~$y_E$
is obtained from~$v$ by a series of operations,  each of which
is like the one described in the above paragraph.
\end{myproof}

In fact,  we have
\begin{proposition}\mylabel{p:oldomination}
$y_E\geq w_E$, where $w_E$
is the element of~$I(d)$ attached as in~\cite[\S2.2.2]{ru} to~$E$.
\end{proposition}
\begin{myproof}
The strategy is similar to that of the proof of 
Proposition~\ref{p:newdomination}.     There corresponds to~$y_E$
(\cite[Proposition~4.3]{kr})
a subset~$\mon_{y_E}$ of~$\andposv$ that is 
`distinguished' in the sense of~\cite[\S4]{kr}.  (Furthermore,
the subset is symmetric about the diagonal and contains evenly
many diagonal elements~\cite[Proposition~5.2.1]{ru}.)

We first give an explicit description of~$\mon_{y_E}$.   Let the elements
of $\projeven{E}$ arranged in decreasing order be
\[ (r_1,r_1^*), \ldots, (r_u,r_u^*), \ldots, (r_t,r_t^*)  \]
where $u$ is such that $(r_u,r_u^*)$ but not $(r_{u+1},r_{u+1}^*)$ belongs
to~$\andposv$, or, equivalently, $r_u>r_u^*$ but $r_{u+1}<r_{u+1}^*$. 
Throughout this proof,  we use~$i$ and $j$ consistently to denote
integers in the range $1$, \ldots, $u$ and $u+1$, \ldots, $t$
respectively.

Clearly $(r_j,r_j^*)$ 
are all horizontal projections.  Let $p(j)$ be such that $(r_j,r_{p(j)}^*)$
belongs to~$E$:    all the column indices of elements of~$E$ must appear
as column indices also in~$\projeven{E}$,  for no vertical projection
is left out in~$\projeven{E}$.   Then $(r_{u+1},r_{p(u+1)}^*)>\ldots
>(r_t,r_{p(t)}^*)$ is a \vchain\ and $p(u+1)<\ldots<p(t)$.

Let $\sigma$ denote the function~$\{u+1,\ldots,t\}\to\{1,\ldots,u\}$ 
defined inductively as follows:   
\begin{itemize}
\item $\sigma(t)$ is largest possible such that $r_t>r_{\sigma(t)}^*$;
\item $\sigma(t-1)$ is largest possible in $\{1,\ldots,t\}\setminus
\{\sigma(t)\}$ such that $r_{t-1}>r_{\sigma(t-1)}^*$;
\item[\vdots] 
\item $\sigma(j)$ is largest possible in $\{1,\ldots,t\}\setminus
\{\sigma(t),\sigma(t-1),\ldots,\sigma(j+1)\}$ such that 
$r_j>r_{\sigma(j)}^*$.
\end{itemize}
Such a choice of~$\sigma$ is possible.  Indeed, 
\begin{enumerate}
\item $\sigma(t)\geq p(t)$, \ldots, $\sigma(j)\geq p(j)$, \ldots,
$\sigma(u+1)\geq p(u+1)$;
\item If $\sigma(j)>p(j)$,  then $\sigma(j-1)\geq p(j)$ (for
$r_{j-1}>r_j>r_{p(j)}^*$).
\end{enumerate}
We have
\begin{multline*}
\mon_{y_E}=\{(r_j,r_{\sigma(j)}^*), (r_{\sigma(j)}, r_j^*)\st u+1\leq j\leq t\}
\bigcup\\
 \{(r_i,r_i^*)\st 1\leq i\leq u, 
\not\exists~j~\textrm{with}~i=\sigma(j)\}
\end{multline*}

Next we draw some conclusions from the above description of~$\mon_{y_E}$:
\begin{enumerate}
\item[(a)]  If $E_1>\ldots>E_\ell$ be the decomposition of~$E$ into
intertwined components,  then $\mon_{y_E}=\projeven{E_1}\cup\cdots\cup
\projeven{E_{\ell-1}}\cup\mon_{y_{E_\ell}}$.
\item[(b)] Vertical projections of all elements preceding the
critical element belong to~$\mon_{y_E}$.
\item[(c)] If there exists an element $\alpha$ in~$E_\ell$ of type~H
(there is at most one such element) and $p_h(\alpha)$ belongs to
$\projeven{E}$,  then $p_h(\alpha)\in\mon_{y_{E_\ell}}$.
\item[(d)] For each $\alpha$ in~$E$ there exists a unique element~$\beta$
in~$\mon_{y_E}$  that shares its column index with~$\alpha$.   This element
lies on or above the diagonal and its row index is no smaller than that
of~$\alpha$.   If~$E$ is~$\alpha_1>\alpha_2>\ldots$,  then the
corresponding elements form a \vchain~$\beta_1>\beta_2>\ldots$ 
in~$\mon_{y_E}$.
\item[(e)] Suppose that $\alpha$ is the critical element of~$E$ 
and $\beta\neq p_v(\alpha)$ where $\beta$ is the corresponding
element in~$\mon_{y_E}$ (see~(d)).    Then $p(j)=\sigma(j)$ $\forall$~$j$ and
$\proj{E} = \projeven{E}$.   
\item[(f)]  Let~$\alpha$ be the critical element of~$E$.
If $\alpha$ has type~V, its horizontal projection $p_h(\alpha)$
belongs to $\projeven{E}$ 
(in other words~$p_h(\alpha)=(r_{u+1},r_{u+1}^*)$),
and $\sigma(j)=p(j)$ $\forall$~$j$, then
the only elements of~$\projeven{E}\cap\andposv$
smaller than~$p_v(\alpha)$ are the vertical projections
of elements of~$E$ (evidently of those beyond the critical element).
\end{enumerate}
\begin{myproof}
(a)  Observe that the critical element~$(r_{u+1},r_{p(u+1)}^*)$ 
belongs to~$E_\ell$
(for the critical element is intertwined with all its successors).  
Since $\sigma(j)\geq p(j)$ for all $j$ and $p(u+1)<\ldots<p(t)$,  the
conclusion follows.

(b) This is because $\{\sigma(t),\ldots,\sigma(u+1)\}\subseteq
\{p(u+1),p(u+1)+1,\ldots,t\}$.    

(c) Let $p_h(\alpha)=(r_s,r_s^*)$.   Since $\alpha$ is not connected to
(but is intertwined with) any of its successors,  we have $r_j\not>r_s^*$
$\forall$ $j$,  so $s\not\in\{\sigma(u+1),\ldots,\sigma(t)\}$.
And clearly $s\leq u$,  so the conclusion follows.

(d)~%
Since $p_v(\alpha)\in\projeven{E}$,  the existence and uniqueness
of~$\beta$ is clear from the description of~$\mon_{y_E}$ above.
Also clear from the description is that the only elements below the
diagonal in~$\mon_{y_E}$ are those with column indices $r_j^*$,  but
$p_v(\alpha)=(r_i,r_i^*)$ for some~$i$ ($p_v(\alpha)\in\andposv$ surely),
so~$\beta$ lies on or above the diagonal.   

To see that the row index of~$\beta$ is no smaller than that of~$\alpha$,
first note that this is clear if $\beta=p_v(\alpha)$.   
If~$\alpha$ precedes the
critical element,   then $\beta=p_v(\alpha)$ by~(b).   
So suppose
that $\alpha=(r_j,r_{p(j)}^*)$ and further that $p(j)=\sigma(j')$ for some
$j'$, $u+1\leq j'\leq t$ 
(if no such $j'$ exists,  then again $\beta=p_v(\alpha)$
by the description of~$\mon_{y_E}$).     
Then $p(j)\geq p(j')$ (for $\sigma(j')\geq
p(j')$),  so $j\geq j'$ (for $p(u+1)<\ldots<p(t)$).   Since $\beta
=(r_{j'},r_{\sigma(j')}^*)$, it follows that $r_{j'}\geq r_j$, i.e.,
$\beta$~has no smaller row index than that of~$\alpha$.   

Finally, that $\beta_1$, $\beta_2$, \ldots
form a \vchain\ follows readily by combining the assertion just
proved with the distinguishedness of~$\mon_{y_E}$.

(e)~The assumption that $\beta\neq p_v(\alpha)$ implies that
$p_v(\alpha)(=(r_{p(u+1)},r_{p(u+1)}^*))$ does not belong to~$\mon_{y_E}$,
which means $p(u+1)=\sigma(j)$ for some~$j$.   If $j>u+1$, we
have $\sigma(j)\geq p(j)>p(u+1)$ (see~(1) above), a contradiction,
so $p(u+1)=\sigma(u+1)$.   By~(2) above,  it follows 
that $p(j)=\sigma(j)$ for all~$j$.

Suppose that $\proj{E}$  has oddly many elements.  Let~$i$
be such that~$(r_i,r_i^*)$ is the vertical projection of the last element,
say~$\lambda$, of~$E$.   Since $p_h(\lambda)\not\in\projeven{E}$,  it follows
that $i>p(t)$ (note that $(r_{p(t)}, r_{p(t)}^*)$ is the vertical projection
of the element of~$E$ with horizontal projection~$(r_t,r_t^*)$).   Since
$r_t>r>r_i^*$,  where $r$ denotes the row index of~$\lambda$,  
we have $\sigma(t)\geq i>p(t)$ contradicting the previous assertion.

(f)~Note that $(r_{p(u+1)}, r^*_{p(u+1)})$ is the vertical projection
of~$\alpha$ (by the definition of~$p$).   
Suppose that there exists~$(r_i,r_i^*)$ with $i> p(u+1)$ that is not
the vertical projection of any element of~$E$, i.e., there does not exist~$j$
with~$i=p(j)$.    Then~$(r_i,r_i^*)$ is a horizontal projection,
evidently of some predecessor of~$\alpha$.   If $r_{u+1}<r_i^*$,  then~$\alpha$
is not connected with that predecessor,  therefore neither to its
immediate predecessor,  and so of type~S (rather than~V as assumed).
We may therefore assume that $r_{u+1}>r_i^*$.   Now, if $i=\sigma(j)$
for some $j>u+1$,  then $\sigma(j)\neq p(j)$, a contradiction;  
if not, then it follows from the definition of~$\sigma$ that
$\sigma(u+1)\geq i>p(u+1)$, again a contradiction.    (It is easy
to construct counter-examples to the assertion with the critical
element being the last element of~$E$ and its horizontal projection
being not in~$\projeven{E}$,  in which case the hypothesis that 
$\sigma(j)=p(j)$ for all $j$ is vacuously satisfied.)
\end{myproof}

We are finally ready for the proof of the proposition.
By~\cite[Lemmas~4.5,~5.5]{kr},  it is enough to show that every \vchain\
in $\mon_E$ is dominated by one in~$\mon_{y_E}$.
Let $E_1>\ldots>E_\ell$ be the decomposition of~$E$ into
intertwined components.
Take a \vchain~$C$ in~$\mon_E$.   As observed in the proof of
Proposition~\ref{p:newdomination},  $C$ is just a concatenation of
\vchains~$C_1$, \ldots, $C_\ell$ with $C_j$ being a \vchain\ in~$\mon_{E_j}$.
We have already seen in Lemma~\ref{l:p:newdomination:1} that there exist
\vchains~$D_1$, \ldots, $D_{\ell-1}$ in $\projeven{E_1}$, \ldots,
$\projeven{E_{\ell-1}}$ respectively dominating $C_1$, \ldots, $C_{\ell-1}$.
In the light of~(a) above,  we'd be done if we can find~$D_\ell$ in
$\mon_{y_{E_\ell}}$  dominating~$C_\ell$,  for then the concatenation
$D_1>\ldots>D_{\ell-1}>D_\ell$ would be a \vchain\ in~$\mon_{y_E}$ dominating~$C$.
As in the proof of Lemma~\ref{l:p:newdomination:2},   we may reduce
to the case when~$C_\ell$ lies weakly above the diagonal (this follows
from the proof of \cite[Lemma~5.15]{gr} and the symmetry about the diagonal
of monomials attached to~\vchains).   

We now show that such a chain~$D_\ell$ exists.  In fact, let us show:
for an intertwined \vchain~$F$ and $\mu_1>\mu_2>\ldots$ a maximal
\vchain\ in~$\mon_F$ lying weakly above the diagonal,  there exists
$\nu_1>\nu_2>\ldots$ in~$\mon_{y_F}$ lying weakly above the diagonal
that dominates $\mu_1>\mu_2>\ldots$.   The goal being
analogous to that of Lemma~\ref{l:p:newdomination:2},  we adopt the
notation and arguments from the first of the three
parts of that proof.   There are two possibilities for 
$\mu_1>\mu_2>\ldots$, namely (*) and (**) as in the 
proof of that lemma.

First consider~(**).   If $p_h(\alpha_k)$ belongs to $\projeven{F}$,
then~(**) is contained in~$\mon_{y_F}$ by~(b) and~(c) above.
If not, then $\alpha_k$ is the last element of~$F$,  so that all
projections of~$F$ belong to~$\andposv$.     In this case, $\mon_{y_F}
=\projeven{F}=\mon_{\spnew{F}}$,  and we're done by invoking
Lemma~\ref{l:p:newdomination:1}.

Now consider the \vchain~(*).    Because of~(b) and~(d) above,
it follows that the \vchain~$\beta_1>\beta_2>\ldots$ as in~(d) 
dominates~(*) except in the following situation:   the critical
element~$\alpha_{k+1}$ has type~V and $\beta_{k+1}\neq p_v(\alpha_{k+1})$.
So assume that we are in this situation
(which means that the types of elements of~$F$ are as in Case~II 
on page~\pageref{case2}
and in particular that $k$ is odd).   Assertions~(e) and~(f) above apply.

The elements $p_v(\alpha_1)$, \ldots, $p_v(\alpha_k)$ belong to $\mon_{y_F}$
(by (b)).   If there is one other element in $\mon_{y_F}$ that dominates
$p_v(\alpha_{k+1})$,  then these elements together form a 
\vchain~$\gamma_1>\ldots>\gamma_{k+1}$ in~$\mon_{y_F}$
that dominates $p_v(\alpha_1)>\ldots>p_v(\alpha_k)>p_v(\alpha_{k+1})$,
and $\gamma_1>\ldots>\gamma_{k+1}>\beta_{k+2}>\beta_{k+3}>\ldots$ dominates
(*),  and we're done.     So assume that this is not the case.
From~(e) and~(f) above it follows that $\proj{F}$ consists precisely of
$p_v(\alpha_1)$, \ldots, $p_v(\alpha_{k})$  and both projections
of $\alpha_{k+1}$, $\alpha_{k+2}$, \ldots,  
and so of an odd number (because $k$ is odd), contradicting~(e).
\end{myproof}

\mysection{Pfaffians and their Laplace-like expansions}\mylabel{s:pfaffian}
This section can be read independently of the rest of the paper.
We define here the {\em Pfaffian\/} of a matrix of even size that
is skew-symmetric along the anti-diagonal and show that it satisfies
a Laplace-like expansion formula similar to the one for the determinant.
In fact we define the Pfaffian by such a formula: see Eq.~(\ref{e.pfaffian}).
We then show that it is independent of the choice of the integer involved
in the expansion and that it is a square root of the determinant 
(Corollary~\ref{c:p:pfaffian}).
The expansion formula is used crucially in the proof of the main
Lemma~\ref{l:main} in~\S\ref{s:proof}. 

\mysubsection{The Pfaffian defined by a Laplace-like expansion}%
\mylabel{ss:pfaffdef}
Let $n$ be a non-negative integer.   For $k$ an integer,  define
$k^*=2n+1-k$.
Let $A=(a_{ij})$ be a $2n\times 2n$ matrix that is
skew-symmetric along the anti-diagonal, meaning that $a_{ij}=-a_{j^*i^*}$
for $1\leq i,j\leq 2n$.     We will be considering submatrices of~$A$.     
Let $A_{r,c}$ denote the submatrix obtained by deleting the row numbered~$r$ 
and the column numbered~$c$; 
$A_{r_1r_2,c_1c_2}$ the submatrix obtained by deleting rows 
numbered $r_1$,~$r_2$ and column numbers $c_1$,~$c_2$; and so forth. 
Let $D$, $D_{r,c}$, $D_{r_1r_2,c_1c_2}$, \ldots denote respectively 
the determinants of~$A$, $A_{r,c}$, $A_{r_1r_2,c_1c_2}$, \ldots.

We define the {\em Pfaffian\/}~$Q$ of the matrix~$A$ by induction on~$n$:
for $n=0$, set $Q:=1$;  for $n\geq 1$, set
\begin{equation}\label{e.pfaffian} 
  Q:=\sum_{j=1}^{2n} (-1)^{m+j^*}
 \sign(mj)\,a_{m,j^*}\, Q_{mj,j^*m^*} \end{equation}
where $m$ is a fixed integer, $1\leq m\leq 2n$;
$Q_{mj,j^*m^*}$ is the Pfaffian of the submatrix $A_{mj,j^*m^*}$; and,
for natural numbers $i$ and $j$,  
\[ \sign(ij):=\left\{
        \begin{array}{cl}
          1 & \textrm{if $i<j$}\\
          -1 & \textrm{if $i>j$}\\
            0 & \textrm{if $i=j$}
          \end{array}\right. \]  
($Q_{mj,j^*m^*}$ is not defined when $j=m$ but this does not matter
since $\sign(mj)=0$ then).
To see that the expression~(\ref{e.pfaffian}) is independent of the choice
of~$m$,  proceed by induction on~$n$.   If $p$ is another
choice,  then, by the induction hypothesis,
$Q_{mj,j^*m^*}$ equals 
\[ \sum_{k=1}^{2n} 
 (-1)^{p+k^*}\sign(pm)\sign(pj)\sign(k^*j^*)\sign(k^*m^*)
 \sign(pk)\,a_{p,k^*} Q_{pmjk,k^*j^*m^*p^*}  \]
and, similarly, 
$Q_{pk,k^*p^*}$ equals
\[ \sum_{j=1}^{2n} 
 (-1)^{m+j^*}\sign(mj)\sign(pm)\sign(mk)\sign(k^*j^*)\sign(j^*p^*)
 \,a_{m,j^*} Q_{pmjk,k^*j^*m^*p^*}  \]
so that, irrespective of whether~$m$ or~$p$ is chosen, we get
\[\begin{array}{r} Q =
\sum_{j, k=1}^{2n} 
 (-1)^{m+j^*+p+k^*}\sign(mj)\sign(pm)\sign(pj)\sign(k^*j^*)\sign(k^*m^*)
 \cdot\\ \sign(pk) \, a_{m,j^*}\, a_{p,k^*}\, Q_{pmjk,k^*j^*m^*p^*}. \end{array} \]

Since \[(-1)^{m+j^*}\sign(mj)a_{m,j^*}Q_{mj,j^*m^*}\] is symmetric in~$m$ 
and~$j$ (for we have
$(-1)^m=-(-1)^{m^*}$, $(-1)^{j^*}=-(-1)^j$, $\sign(mj)=-\sign(jm)$,
$a_{m,j^*}=-a_{j,m^*}$, and, obviously, $Q_{mj,j^*m^*}=Q_{jm,m^*j^*}$), the
summation in equation~(\ref{e.pfaffian}) can be taken over $m$:
\begin{equation}\label{e.pfaffian.col} 
  Q=\sum_{m=1}^{2n} (-1)^{m+j^*}
 \sign(mj)\,a_{m,j^*}\, Q_{mj,j^*m^*} \end{equation}
\begin{corollary}\mylabel{c:pfterms}
The number of terms in the Pfaffian of a generic $2n\times 2n$ matrix
skew-symmetric along the anti-diagonal is $(2n-1)\cdot(2n-3)
\cdot\cdots\cdot3\cdot1$.
\textup{By convention we take this number to be~$1$ when $n=0$ (in analogy with
  the convention $0!=1$).}\hfill$\Box$
  \end{corollary}
\mysubsection{Pfaffians and determinants}\mylabel{ss:pfaffdet}
\begin{proposition}\mylabel{p.pfaffian}
For integers $a$, $j$, $k$ such that $1\leq a, j, k\leq 2n$ and
$a\neq j$, $a\neq k$, 
\[ D_{aj,k^*a^*}=(-1)^{n-1}Q_{aj,j^*a^*}Q_{ak,k^*a^*}. \]
\end{proposition}
\begin{myproof}
Proceed by induction.  
Writing the Laplace expansion for~$D_{aj,k^*a^*}$ along row $k$ of $A_{aj,k^*j^*}$,
we get
\[ D_{aj,k^*a^*}= \sum_{i=1}^{2n} (-1)^{k+i^*}\sign(ak)\sign(jk)\sign(i^*k^*)\sign(i^*a^*)\, a_{k,i^*} D_{ajk,i^*k^*a^*}. \]
Writing the Laplace expansion for $D_{ajk,i^*k^*a^*}$ along column $j^*$ of
$A_{ajk,i^*k^*a^*}$,  we get
\[\begin{array}{r} D_{ajk,i^*k^*a^*}= \sum_{\ell=1}^{2n} (-1)^{\ell+j^*}
\sign(a\ell)\sign(j\ell)\sign(k\ell)
\sign(i^*j^*)\sign(k^*j^*)\cdot \\ \sign(j^*a^*)
\, a_{\ell,j^*} D_{ajk\ell,i^*k^*j^*a^*}.\end{array} \]
By the induction hypothesis,
\[ D_{ajk\ell,i^*k^*j^*a^*}=(-1)^{n-2}Q_{ajk\ell,\ell^*k^*j^*a^*}Q_{ajki,i^*k^*j^*a^*}\]
Substituting this into the expression for $D_{ajk,i^*k^*a^*}$ and
the result in turn into the expression for $D_{aj,k^*a^*}$,
and rearranging terms---we have replaced $\sign(i^*k^*)$ by $\sign(ki)$ and
$(-1)^{n-2}\sign(jl)$ by $(-1)^{n-1}\sign(lj)$---we get 
\[\begin{array}{l}
D_{aj,k^*a^*} =  (-1)^{n-1}\cdot\\
\left(
\sum_{i=1}^{2n}((-1)^{k+i^*}\sign(ak)\sign(jk)\sign(i^*j^*)\sign(i^*a^*))\, \sign(ki) a_{k,i^*} Q_{ajki,i^*k^*j^*a^*}
                      \right)\\
 \left(
\sum_{\ell=1}^{2n}((-1)^{\ell+j^*}\sign(a\ell)\sign(k\ell)\sign(k^*j^*)\sign(j^*a^*)) \sign(\ell j) a_{l,j^*}  Q_{ajk\ell,\ell^*k^*j^*a^*}
           \right)
\end{array}\]
By equations~(\ref{e.pfaffian}) and~(\ref{e.pfaffian.col}),
the factors in the second and third lines of the above display 
are respectively~$Q_{aj,j^*a^*}$ and $Q_{ak,k^*a^*}$, so we are done.
\end{myproof}
\begin{corollary}\mylabel{c:p:pfaffian}
$D=(-1)^nQ^2$.
\end{corollary}
\begin{myproof}
  Put $j=k$ in the proposition.
\end{myproof}

\mysection{The proof}\mylabel{s:proof}\mylabel{s:lemma}
We are now ready to prove our result (Theorem~\ref{t:main}).
Lemma~\ref{l:main} is the technical result that enables the proof.
Its proof uses the results of~\S\ref{s:newforms},~\ref{s:pfaffian}.

Notation is fixed as in~\S\ref{ss:theorem}.   
\mysubsection{Setting it up}\mylabel{ss:setupproof}
Our goal is to prove: 
\begin{quote}\begin{sl}
  Every monomial in~$\rootsv$ that is not \ortho-dominated by~$w$
occurs as an initial term with respect to the term order~$\torder$ of
an element of the ideal~$I$ of the tangent cone.
\end{sl}\end{quote}
As explained in~\S\ref{ss:theorem},  putting this assertion
together with the main result of~\cite{ru} yields Theorem~\ref{t:main}.

Let~$I'$ be the ideal generated by $f_\tau$, 
$\tau\in I(d)$, $v\leq \tau\not\leq w$.    Since~$I'\subseteq I$, and since
a monomial in~$\rootsv$ that is not \ortho-dominated by~$w$ contains, by the
definition of \ortho-domination (\S\ref{ss:dominate}),  a \vchain\ in~$\posv$
that is not \ortho-dominated by~$w$, 
it suffices to prove the following (after which it will follow
that $I'=I$):
\begin{quote}\begin{sl}
  Every \vchain\ that is not \ortho-dominated by~$w$ 
occurs as the initial term of an element of~$I'$.
\end{sl}
\end{quote}
Putting $j=1$ in~Lemma~\ref{l:main} below yields this,
so it suffices to prove that lemma.
\mysubsection{The main lemma}\mylabel{ss:lemma}

\newcommand\chain{A}
Fix a \vchain~$\chain:\alpha_1>\ldots>\alpha_m$ that is not \ortho-dominated
by~$w$.    Let $j$ be an integer, $1\leq j\leq m$.   Define~$A_j$ to be
the sub-\vchain\ $\alpha_1>\ldots>\alpha_j$.  Set
\[\Gamma_j:=\left\{
\begin{array}{cl}
\projeven{\chain_{j}} & \textup{if $\card\proj{\chain_{j}}$ is odd}\\
\projeven{\chain_{j}}\setminus\{p_v(\alpha_{j}),p_h(\alpha_j)\} & \textup{if
 $\card\proj{\chain_j}$ is even}\\
\end{array}
\right.\]
See~\S\ref{ss:nf:prepare} for the definition of~$\proj{}$ and $\projeven{}$.
Observe that 
\begin{quote}
(\ddag)\quad
if~$\card\proj{\chain_{j-1}}$ is even (equivalently $\proj{\chain_{j-1}}=\projeven{\chain_{j-1}}$),  then 
$\Gamma_j=\projeven{\chain_{j-1}}$, no matter whether $\card\proj{\chain_j}$
is even or odd.
\end{quote}
$\Gamma_j$ being a subset of even cardinality, say $2q_j$,
 of the diagonal elements
of~$\rootsv$,  it defines an element of~$I(d)$.   The corresponding
Pfaffian we denote by~$f_j$.   The degree of~$f_j$ is
$q_j$ and the number of terms in~$f_j$ is, by Corollary~\ref{c:pfterms},
$n_j:=(2q_j-1)\cdot(2q_j-3)\cdot\cdots\cdot3\cdot1$.    By convention, $n_j=1$
when $q_j=0$.
\begin{lemma}\mylabel{l:main}
Let~$A: \alpha_1>\ldots>\alpha_m$ be a \vchain\ not \ortho-dominated by~$w$.
For every integer~$j$, $1\leq j\leq m$,  there exists a homogeneous 
element~$F_j$ of the ideal~$I'$ 
such that 
\begin{enumerate}
\item  For a monomial occurring with non-zero coefficient in~$F_j$, consider
the set (counted with multiplicities) of the projections on the diagonal 
of the elements of~$\rootsv$ that occur in the monomial.   This set 
is the same for every such monomial.
\item
The sum of the initial $n_j$ terms (with respect to
the term order~$\torder$) of~$F_j$ is $f_jX_{\alpha_j}\cdots X_{\alpha_m}$.  
\end{enumerate}
Consider any fixed monomial (occurring with non-zero coefficient) in~$F_j$
other than one in $f_jX_{\alpha_j}\cdots X_{\alpha_m}$.  From (1) and (2) 
it follows that,  given an integer~$b$,  $j\leq b\leq m$,   
there exists precisely one $X_{\delta_b}$ occurring
in the monomial with the row index of~$\delta_b$ being that of~$\alpha_b$.
\begin{enumerate}
\item[3]
There exists~$b$ for which $\delta_b\neq \alpha_b$ and, for
the largest~$b$ of this kind, either $\delta_b\not\in\posv$ or
the column index of~$\delta_b$ is less than
that of~$\alpha_b$.
\end{enumerate}
\end{lemma}
\begin{myproof}
  Proceed by an induction on~$m$ and then another (in reverse) on~$j$.
Let us suppose that we know the result for~$j$ and prove it for~$j-1$.
The proof below covers also the base cases for the induction.
Consider~$\proj{\chain_{j-1}}$.   

Suppose first that its cardinality
$\card\proj{\chain_{j-1}}$ is odd.
Write~$\chain$ as $C>D$ with $C=\chain_{j-1}$ and $D$ being 
$\alpha_j>\ldots>\alpha_m$.
Observe that the last intertwined
component of~$C$ has at least two elements.    Let $\new{\chain}$ be the
new form of~$A$ constructed as in~\S\ref{ss:nf:const}.   
Since~$\new{A}$ has fewer elements than~$A$ (Proposition~\ref{p:fewer})
and is not \ortho-dominated by~$w$ (Proposition~\ref{p:newdomination}),
the induction hypothesis applies to~$\new{A}$.   Apply it with
$k=\card\new{C}+1$ in place of~$j$ in the statement of the lemma.
If~$F$ is the element in~$I'$ as in its conclusion,   set 
$F_{j-1}=X_{\alpha_{j-1}}F$.

We claim that $F_{j-1}$ has the desired properties.
That it satisfies~(1) is clear.   We now observe that it satisfies~(2).
Since $\proj{{\new{\chain}}_{k-1}}=\proj{\new{C}}$ has evenly many
elements (Proposition~\ref{p:fewer}),   it follows (observation~(\ddag) above)
that $\Gamma_k$ (calculated for $\new{A}:\new{C}>D$)  equals 
$\projeven{\new{C}}=\proj{\new{C}}$.
On the other hand, $\Gamma_{j-1}=\projeven{C}=\proj{\new{C}}$ (since
$\proj{\chain_{j-1}}$ is odd,  by Proposition~\ref{p:fewer}).
So~$F_{j-1}$ satisfies~(2).   That $F_j$ satisfies~(3) is readily verified.

Now suppose that $\card\proj{\chain_{j-1}}$ is even.  
Apply the induction hypothesis with~$j$ and let~$F_j$ be as in its
conclusion. 
The base case $j-1=m$ needs to be treated separately here, as follows.
Let~$y_A$ be the element of~$I(d)$ defined as in~\S\ref{ss:y}.
We take $F_j$ to be the Pfaffian~$f_{y_A}$ attached to~$y_A$
(see~\S\ref{ss:tangentcone}). 
That $F_j$ belongs to~$I'$ follows from Propositions~\ref{p:ygeqv}
and~\ref{p:oldomination}.    The rest of the proof is the same for the
induction step as well as the base case.

From the observation (\ddag) above, it follows that $\Gamma_j=
\proj{\chain_{j-1}}$.     Here is a picture of~$\Gamma_j$ (the
solid circles denote elements of~$\Gamma_j$):
\setlength{\unitlength}{.6cm}
\\[1mm]
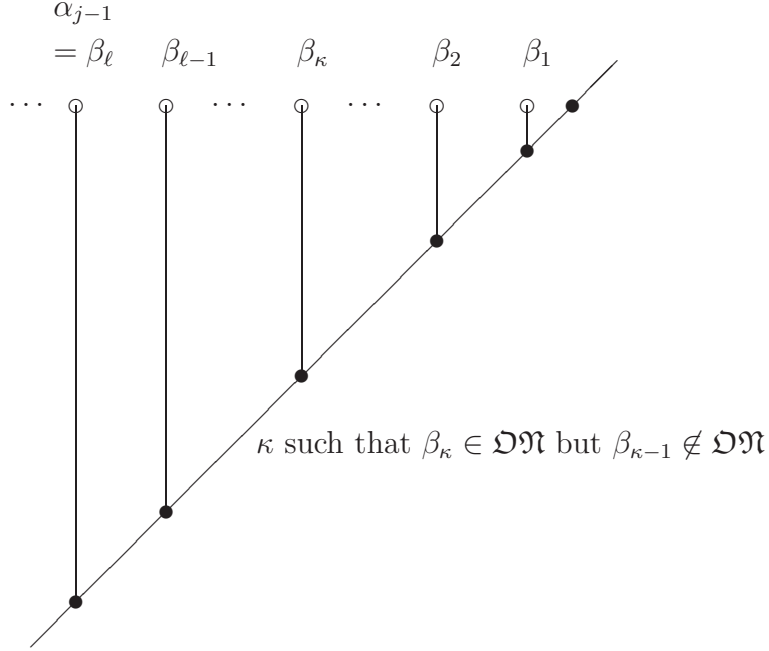
\begin{picture}(18,14)(-6,-1)
  \put(-2,-2){\line(1,1){13}}
\put(-1,-1){\circle*{.3}}
\put(1,1){\circle*{.3}}
\put(4,4){\circle*{.3}}
\put(7,7){\circle*{.3}}
\put(9,9){\circle*{.3}}
\put(10,10){\circle*{.3}}

\put(0.9,11){$\beta_{\ell-1}$}
\put(1,10){\circle{.3}}
\put(-1,10){\circle{.3}}
\put(2,10){\ldots}
\put(4,10){\circle{.3}}
\put(5,10){\ldots}
\put(7,10){\circle{.3}}
\put(9,10){\circle{.3}}

\linethickness{.1mm}
\put(-1,-1){\line(0,1){11}}
\put(1,1){\line(0,1){9}}
\put(4,4){\line(0,1){6}}
\put(7,7){\line(0,1){3}}
\put(9,9){\line(0,1){1}}

\put(8.9,11){$\beta_1$}
\put(6.9,11){$\beta_2$}
\put(3.9,11){$\beta_{\kappa}$}
\put(-1.5,11){$=\beta_\ell$}
\put(-1.5,12){$\alpha_{j-1}$}
\put(-2.5,10){\ldots}

\put(3,2.3){$\kappa$ such that $\beta_\kappa\in\posv$ but $\beta_{\kappa-1}\not \in\posv$}
\end{picture}\\[.8cm]
Applying to~$f_j$ the Laplace-like expansion formula (\ref{e.pfaffian}) 
for Pfaffians,    we see that the sum of its initial $n_{j-1}$ terms,
the next $n_{j-1}$ terms, \ldots\ are (up to sign factors) 
$g_\kappa X_{\beta_\kappa}$,
$g_{\kappa+1}X_{\beta_{\kappa+1}}$, \ldots, 
$g_{\ell-1}X_{\beta_{\ell-1}}$, $g_\ell X_{\alpha_{j-1}}$,
\ldots,  where $g_i$ is the Pfaffian associated 
to $\Gamma_{j}\setminus\{p_v(\beta_i),p_h(\beta_i)\}$, so that
the corresponding initial terms of~$F_j$ are 
$g_\kappa X_{\beta_\kappa}X_{\alpha_j}\cdots X_{\alpha_m}$,
$g_{\kappa+1}X_{\beta_{\kappa+1}}X_{\alpha_j}\cdots X_{\alpha_m}$, \ldots, 
$g_{\ell-1}X_{\beta_{\ell-1}}X_{\alpha_j}\cdots X_{\alpha_m}$, $g_\ell X_{\alpha_{j-1}}X_{\alpha_j}\cdots X_{\alpha_m}$,
\ldots.   We will 
now modify~$F_j$ (by subtracting from it elements of~$I'$)
so as to kill the terms $g_\kappa X_{\beta_\kappa}X_{\alpha_j}\cdots X_{\alpha_m}$, 
\ldots,
$g_{\ell-1}X_{\beta_{\ell-1}}X_{\alpha_j}\cdots X_{\alpha_m}$.   
But of course this needs to be done carefully in order that the resulting
element of~$I'$ has the desired properties.

Write $\chain$ as $C>D$ where $C=\chain_{j-1}$ and $D$ 
is $\alpha_j>\ldots>\alpha_m$.
We may assume that the last intertwined component of~$C$ consists
of at least two elements,  for otherwise
~$F_j$ itself without further modification has the
desired properties (we can take~$F_{j-1}$ to be~$F_j$).   We may further assume
that there is some element of $\proj{A_{j-1}}$ that is strictly in between the
vertical and horizontal projections of~$\alpha_{j-1}$, for otherwise again
we can take~$F_{j-1}$ to be~$F_j$.
Consider the new forms of~$\chain$ as in~\S\ref{ss:nf:const}. 
In their construction there is the choice involved of a diagonal element
strictly in between the vertical and horizontal projections of the last
element of $C$.    We can choose this element to be the vertical projection
of~$\beta_i$ where $\kappa\leq i\leq \ell-1$.     Corresponding to each choice
we get a new form which let us denote $\new{\chain}(i)$ ($=\new{C}(i)>D$).
Since~$\new{A}(i)$ has fewer elements than~$A$ (Proposition~\ref{p:fewer})
and is not \ortho-dominated by~$w$ (Proposition~\ref{p:newdomination}),
the induction hypothesis applies to~$\new{A}(i)$.
Apply it with 
$k=\card\new{C}(i)+1$ in place of~$j$ in the statement of the lemma.
Let $F(i)$ in~$I'$ be as in its conclusion.   
Set~$F_{j-1}=F_{j}-\sum_{i=\kappa}^{\ell-1} F(i)X_{\beta_i}$.

It remains only to verify that~$F_{j-1}$ has the desired properties.
Since $\proj{{\new{\chain}(i)}_{k-1}}=\proj{\new{C}(i)}$ has evenly many
elements (Proposition~\ref{p:fewer}),   it follows (observation~(\ddag) above)
that $\Gamma_k$ (calculated for $\new{A}(i):\new{C}(i)>D$)  equals 
$\projeven{\new{C}(i)}=\proj{\new{C}(i)}$.
From the definition of~$\new{C}(i)$ and observation~(\ddag), it
follows that $\proj{\new{C}(i)}$
is $\Gamma_j\setminus\{p_v(\beta_i),p_h(\beta_i)\}$.
So the sum of the initial~$n_{j-1}$ terms of~$F(i)$ is~$g_iX_{\alpha_j}\cdots X_{\alpha_m}$.
That~$F_{j-1}$ has the desired properties can now be readily verified.
\end{myproof}

	\newcommand\citenumfont[1]{\textbf{#1}}

\bibliographystyle{bibsty-final-no-issn-isbn}
\addcontentsline{toc}{section}{References}
\ifthenelse{\equal{\finalized}{no}}{
\bibliography{abbrev,references}
}{
}

\end{document}